\date{} % % March 17 2003
\title{Exploration trees and conformal loop ensembles}
\author{Scott Sheffield \thanks{Courant Institute, NYU.  Partially supported
by NSF grant DMS 0403182.}}
\documentclass[11pt,naturalnames]{article}
\usepackage{amsmath,amssymb,amsthm,amsfonts,graphicx,
xspace,url,picinpar,graphicx}

\input labelfig.tex
\input epsf.tex

\newif\ifhyper\IfFileExists{hyperref.sty}{\hypertrue}{\hyperfalse}

\ifhyper\usepackage{hyperref}\fi

\newif\ifdraft
\drafttrue
\def\note#1/{\ifdraft {\bf [#1]}\fi}
\long\def\comment#1{}
\numberwithin{equation}{section}
\numberwithin{figure}{section}

\newtheorem{theorem}{Theorem}
\numberwithin{theorem}{section}

\newtheorem{lemma}[theorem]{Lemma}
\newtheorem{proposition}[theorem]{Proposition}
\newtheorem{conjecture}[theorem]{Conjecture}
\newtheorem{problem}[theorem]{Problem}

\theoremstyle{remark}
\theoremstyle{remark}

\newcommand{\R}{\mathbb{R}}
\newcommand{\C}{\mathbb{C}}

\def\H{\mathbb{H}}
\def\D{\mathbb{D}}

\def\L{\mathcal{L}}
\def\B{\mathcal{B}}
\def\A{\mathcal{A}}

\def\capacity{\mathrm{cap}_\infty}

\def\SLEkk#1/{$\mathrm{SLE}(#1)$}
\def\SLEk/{\SLEkk{\kappa}/}
\def\SLE/{$\mathrm{SLE}$}
\def\SLEr#1/{$\mathrm{SLE(\kappa;#1)}$}
\def\SLEkr/{$\mathrm{SLE(\kappa;\rho)}$}

\def\SLEbkk#1/{$\mathrm{SLE^\mu_\beta}(#1)$}
\def\SLEbk/{\SLEbkk{\kappa}/}
\def\SLEb/{$\mathrm{SLE}^\mu_\beta$}
\def\SLEbr#1/{$\mathrm{SLE^\mu_\beta(\kappa;#1)}$}
\def\SLEbkr/{$\mathrm{SLE^\mu_\beta(\kappa;\rho)}$}

\def\CLEkk#1/{$\mathrm{CLE}(#1)$}
\def\CLEk/{\CLEkk{\kappa}/}
\def\CLEtwo/{\CLEkk2/}
\def\CLE/{$\mathrm{CLE}$}
\def\CLEbkk#1/{$\mathrm{CLE^\mu_\beta}(#1)$}
\def\CLEbk/{\CLEbkk{\kappa}/}

\def\Ito/{It\^o}

\def\proofof#1{{ \medbreak \noindent {\bf Proof of #1.} }}

\def\BESd/{$\mathrm{BES}^\delta$}
\def\BESdx/{$\mathrm{BES}^\delta_x$}
\def\SBESd/{$\mathrm{SBES}^\delta$}
\def\BESQd/{$\mathrm{BESQ}^\delta$}
\def\PV{\mathrm{P.V.}}
\def\sign{\mathrm{sign}}

\def\SmirnovPerc{MR1851632}

\def\RohdeSchramm{MR2153402}
\def\SheffieldWerner{SheffieldWerner}
\def\LSWUST{MR2044671}
\def\WernerLoops{MR2023758}
\def\SkewBesselSurvey{MR2190302}

\def\RevuzYor{MR2000h:60050}

\def\LawlerConformalBook{MR2129588}

\def\CamiaNewman{math.Pr/0504036}
\def\Watanabe{MR1335470}

\def\NienhuisKagerSurvey{MR2065722}

\def\SchrammWilson{math.Pr/0505368}
\def\SchrammSheffieldWilson{SchrammSheffieldWilson}
\def\AizemanBurchard{MR1712629}

\def\PitmanPart{MR1466546}

\begin{document} \maketitle \begin{abstract}

%We exhibit similar correspondences between random continuous analogs of these four objects: the
%{\bf Gaussian loop ensembles} \CLEk/ (introduced here) with $8/3 \leq \kappa \leq 8$, a branching
%variant of \SLEkk \kappa, \kappa -6/ (also introduced here), the {\bf Gaussian free field} (GFF),
%and the tree of ``{\bf flow lines} of the complex vector field $e^{ih/\chi}$,'' where $\chi$ is a
%$\kappa$-dependent constant and $h$ is an instance of the GFF.

We construct and study the conformal loop ensembles \CLEk/, defined for $8/3
\leq \kappa \leq 8$, using branching variants of \SLEk/ called
exploration trees. The \CLEk/ are random collections of countably
many loops in a planar domain that are characterized by
certain conformal invariance and Markov properties.  We conjecture
that they are the scaling limits of various random loop models
from statistical physics, including the $O(n)$ loop models.

%including $O(n)$ loop ensembles, double dimer loops, and Ising
%model cluster boundaries.
%We also define {\bf height functions}, {\bf flow line ensembles},
%and {\bf exploration trees} for these discrete loop ensembles and
%make conjectures about their scaling limits.
%We also introduce \CLEkk \kappa, \rho/, a more general family of loop ensembles.  When $\kappa =
%4$, we interpret these ensembles as ``contour lines of the GFF'' at height spacings that depend on
%$\rho$.
\end{abstract}

\newpage
\tableofcontents
\newpage

\section{Introduction}
\subsection{Overview}
Many two dimensional statistical physical models can be interpreted
as random collections of disjoint, non-self-intersecting loops in a
planar lattice.  For example, the loops may be the boundaries
between the plus spin and minus spin clusters in an Ising model with
spins defined on the faces of a three-regular planar graph.

When the boundary conditions of a random loop model on a simply
connected planar domain are set up so that, in addition to the
loops, there is one chordal path connecting a pair of boundary
points, it is often natural to conjecture (and sometimes
possible to prove) that as the grid size gets finer,
the law of this random path converges to the law of the chordal Schramm-Loewner
evolution \SLEk/ for some $\kappa>0$.
%This is because the \SLEk/ are the only random paths with a certain natural
%conformal Markov property (see Section \ref{s.SLEdefinition}).

Given that such a conjecture holds, it is natural to expect the
collection of all loops to have a scaling limit.  The primary purpose
of this paper is to introduce and study a natural family of candidates for
this scaling limit, called the {\bf conformal loop ensembles} \CLEk/.

The \CLEk/, defined for $8/3 \leq \kappa \leq 8$, are random collections
of loops in a planar domain which look, in some sense, like \SLEk/ locally.
At the extremes, \CLEkk 8/ almost surely consists of a single space-filling loop, which is the
scaling limit of the outer boundary of the free uniform spanning
tree (see \cite{\LSWUST}), and \CLEkk {8/3}/ almost surely contains
no loops at all.  When $8/3 < \kappa < 8$, the collection of loops
in a \CLEk/ is almost surely countably infinite.  When $\kappa = 6$, it is equivalent to the random collection
of loops described in \cite{\CamiaNewman}, where it was
shown to arise as a scaling limit of the cluster boundaries of site
percolation on the triangle lattice.  Like \SLEk/,
the \CLEk/ loops intersect the boundary of the domain almost surely
if and only if $\kappa > 4$.

We will show that if $\L$ is any random collection of loops in the closure of a planar domain $D$ that
satisfies certain natural hypotheses related to conformal invariance
(and if at least one loop in $\L$ intersects $\partial D$ with positive probability)
then $\L$ must be a \CLEk/ for some $4 < \kappa < 8$ (see Theorem \ref{t.cleuniqueness}).
In a separate joint paper with Werner, we will prove
that the \CLEk/ for $8/3 \leq \kappa \leq 4$ are the only random
ensembles of simple loops in $D$ that possess certain (somewhat different)
conformal symmetries (see Section \ref{s.sequels}) \cite{SheffieldWerner}.

A secondary purpose of this paper is to formulate a series of
conjectures and open questions related to conformal loop ensembles.
For example, we will conjecture a continuum analog of the FK cluster
expansion for Potts models and formulate a precise conjecture
about the scaling limit of the $q$-state Potts model for $q \in \{2,3,4 \}$.
We will define height functions for discrete loop ensembles
and conjecture a connection with the Gaussian free field.  We
will conjecture scaling limits for site percolation on certain
random graphs, and we will ask about the continuum fields
that have \CLE/ loops as level lines.

\subsection{Outline}
In Section \ref{s.discretesection}, we will discuss
random discrete collections of loops on planar graphs.
The results in this section (besides the definition of \SLE/) are not
logically necessary for the definition of \CLE/, but they will clarify
our motivation. (The reader who is primarily interested in the continuum
may skim this section on a first reading.)  For simplicity, we will focus on
hexagonal lattice graphs and the so-called $O(n)$ loop
models. We will associate to each disjoint collection of
non-self-intersecting loops a spanning tree of the graph, called an
``exploration tree,'' and show that an appropriate class of trees
is in one-to-one correspondence with the set of disjoint simple loop ensembles.

In Section \ref{s.besselsection} we will assemble some basic facts
about Bessel processes, L\'evy skew stable processes, and \SLEkr/
processes.  In particular, we will argue that radial \SLEr \kappa-6/ and its
variants are the most natural candidates for the limiting laws of a branch of the
exploration tree. We will give a one-to-one correspondence between
the strictly stable L\'evy processes and a family of conformally
invariant continuum exploration path models.

Section \ref{s.besselsection} will also present Conjecture \ref{c.reversible},
which states that when $4 < \kappa < 8$, both chordal \SLEk/ and chordal \SLEr \kappa-6/ are
random continuous paths whose laws are invariant under anticonformal maps
that the swap their endpoints.
This ``time reversal symmetry'' of chordal \SLEk/ processes is necessary to
prove certain symmetries of \CLEk/, discussed below.  The result is known when
$\kappa=6$ (note that \SLEr \kappa-6/ is the same as \SLEk/ in this case), where
the time reversal symmetry follows from the fact that \SLEk/ is a scaling
limit of discrete models that have this symmetry.  We hope that a more general proof
of Conjecture \ref{c.reversible} will appear soon.  In this paper, we will
construct the \CLEk/ and derive some basic properties without using
Conjecture \ref{c.reversible} (except in Section \ref{s.confsymmetry},
where we describe additional facts about the \CLEk/ that could be derived if
Conjecture \ref{c.reversible} were proved).

Section
\ref{CLEsection} will use a coupling of \SLEr \kappa-6/ processes
with different target points to construct a continuum analog of the
exploration tree and use this continuum tree to construct the
\CLEk/.  This approach to defining loops is related to
the one given in \cite{\CamiaNewman} for the case $\kappa=6$, which
also uses variants of \SLEk/ to ``explore'' segments of loops, but it is
somewhat more canonical in that there are fewer arbitrary choices
in the exploration process.

The exploration tree appears to have connections to the Gaussian
free field and to percolation on the so-called discrete gaskets (which we discuss as
open problems in Section \ref{s.openproblems}).  One intriguing point is that
the family of conformally
invariant exploration tree structures will turn out to have a somewhat different
character when $\kappa=4$ and
when $\kappa\not=4$.  This is related to the fact that the family of strictly
stable L\'evy processes corresponding to $\alpha = 1$ has a different character
from the family
corresponding to $\alpha \not = 1$ (see Section \ref{s.besselsection}).

Section \ref{s.confsymmetry} focuses on the non-simple, non-space-filling case $4 < \kappa < 8$
and formulates and proves a uniqueness theorem which says that
any random boundary-intersecting loop ensemble which satisfies certain natural
hypotheses related to conformal invariance
must be a \CLEk/ for some $4 < \kappa < 8$ (Theorem \ref{t.cleuniqueness}).
Conversely, Theorem \ref{t.cleuniqueness} also states that---if Conjecture \ref{c.reversible} is true---
the \CLEk/ themselves satisfy these hypotheses.  (We will not address the analogous questions
for $\kappa \leq 4$ because we expect them to be addressed in a subsequent paper \cite{\SheffieldWerner}.)

Section \ref{s.approx} will describe various approximations of
\SLEkr/ processes and use them to prove Mobius
invariance of \SLEkr/ and other results.  The invariance results are
similar to those in \cite{\SchrammWilson}, but there are technical
issues that arise when the driving parameters of the Loewner
evolutions are not semimartingales; one reasonably simple way around
this involves the approximations mentioned above.  The approximations
also provide the intuition behind some of the conjectures in
Section \ref{s.openproblems}.

Section \ref{s.otherdiscrete} will explore some additional
combinatorial constructions in the discrete setting. In particular,
we use the ``winding number'' of the exploration tree to construct a
height function for each discrete loop ensemble.

Finally, Section \ref{s.openproblems} will present a list of conjectures
and open problems relevant to \CLE/.

\subsection{Planned sequels} \label{s.sequels}

We now mention briefly some work in progress for which we expect
this paper to be a prerequisite. The random closed set $\Gamma$
consisting of points which are not surrounded by a loop in an
instance of \CLEk/ is called the {\bf \CLEk/ gasket}.  The physics
literature contains many non-rigorous calculations about $O(n)$
model scaling limits that are based on conformal invariance
hypotheses (similar to those of Theorem \ref{t.cleuniqueness}, but
not always so explicitly formulated).  It is natural to interpret these results
as predictions about properties of \CLEk/ and \CLEk/ gaskets. In a joint paper with Schramm
and Wilson we will compute the probabilistic fractal dimension of
$\Gamma$ (which agrees with a calculation in the physics literature
made by Duplantier \cite{D:vesicle}) and the distribution of the set of conformal
radii of the loops surrounding a fixed point in the domain (which agrees
with a calculation in the physics literature made by Cardy and
Ziff \cite{CZ})
\cite{\SchrammSheffieldWilson}.  These results will be derived from
Proposition \ref{gasketconfrad}, which we prove here, but which
was first formulated with Schramm and Wilson as part of this joint project.

In a joint paper with Werner we will show that the sets of outermost
loops of the \CLEk/ defined here (i.e., the loops that are not
surrounded by any other loops) for $8/3 \leq \kappa \leq 4$ are the
{\em only} random ensembles $\mathcal L$ of pairwise disjoint,
non-nested simple loops in $D$ with the following natural Markov property:
if $B \subset D$ is a deterministic closed set with simply connected
complement --- and $\tilde B$ is defined to be the closure of the
set of points surrounded by loops that intersect $B$ --- then given
$\tilde B$, the conditional law of the loops in each component of $D
\setminus \tilde B$ is the same as the original law of $\mathcal L$
conformally mapped to that component \cite{\SheffieldWerner}.  We
will also show that the set of outermost loops in a \CLEk/ has the
same law as the set of loop soup cluster boundaries for a loop soup
of intensity $c$ where $c=(3\kappa-8)(6-\kappa)/2\kappa$.  A form of this
statement and a partial proof appear in earlier work by Werner
\cite{\WernerLoops}.

%Work with Schramm will obtain \CLEkk 4/ level lines of the Gaussian
%free field at appropriately spaced intervals.

\medbreak {\noindent\bf Acknowledgments.}
We thank Oded Schramm, Wendelin Werner, and David Wilson, with whom the author
has collaborated on related projects, and without whom this work would not
have been possible.  We also thank Federico Camia, John Cardy, Julien Dub\'edat,
and Charles Newman for many useful conversations.

\section{Discrete motivation} \label{s.discretesection}
\subsection{Exploration trees}  \label{explorationtreesection}
Let $\mathcal H$ be the infinite hexagonal lattice embedded in
$\R^2$. A graph $G$ consisting of the edges and vertices incident to
a finite simply connected subset of the hexagonal faces of $\mathcal
H$ is called a {\bf hexagon graph}. Let $F$, $E$, and $V$ denote,
respectively, the sets of faces, edges, and vertices of $G$.

Let $A$ be an arbitrary subset of $F$; we will refer to the members
of $A$ as {\bf black} hexagons and the members of $F \backslash A$
as {\bf white} hexagons.  Unless otherwise stated, we will also
refer to the hexagons of $\mathcal H$ outside of $F$ as white.  The
components of the boundaries of the black clusters form an ensemble
of non-self-intersecting loops in $G$ which do not intersect one
another.  It is easy to see that this gives a one-to-one
correspondence between subsets $A \subset F$ and collections of
disjoint simple loops in $G$.

Fix a vertex $v_0$ on the outer boundary of $G$ which is incident to
only two edges of $G$, and fix a directed edge $e_0$ outside of $G$,
beginning at some vertex $v_{-1}$ and pointing towards $v_0$. For
each vertex $v$ of $G$, the {\bf exploration path} $T_v(A)$ is a
directed, non-self-intersecting path $v_0, v_1, v_2 \ldots$ that
ends at $v$.  Each $v_k$, for $k \geq 1$, is chosen in such a way that the sequence
$v_{k-2}, v_{k-1}, v_k$ describes a right turn when the directed edge
$(v_{k-2}, v_{k-1})$ points to a black face and a left turn if
$(v_{k-2}, v_{k-1})$ points to a white face {\em unless} this choice
of $v_k$ would fail to lie in the same connected component of $V
\backslash \{v_0, v_1, \ldots, v_{k-1} \}$ as $v$, in which case the
path turns the other direction.

\begin{figure}[t]\label{explorationtreefigure}
%\epsfbox[-45 43 100 240]
\begin{center}
\rotatebox{180}{\reflectbox{\includegraphics{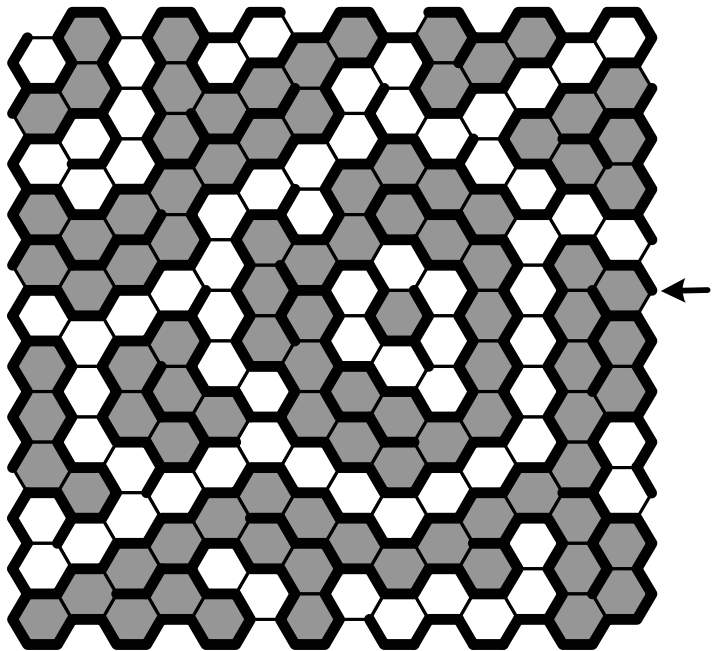}}}
\caption {A coloring of a hexagon graph and the corresponding
exploration tree.  The arrow points to the root vertex $v_0$.}
\end{center}
\end{figure}

The {\bf exploration tree} $T(A)$ of $A$ is the union over all $v\in
V$ of $T_v(A)$. The reader may check that $T(A)$ is in fact an
out-directed spanning tree of $G$, rooted at $v_0$. Readers familiar
with computer algorithms may recognize $T(A)$ as the depth-first
search tree of $G$ beginning at $v_0$ under the rule that one
searches right first after tracing an edge directed towards a black face and
left first after tracing an edge directed towards a white face. See Figure
\ref{explorationtreefigure}.

A vertex $v$ is called a {\bf branch point} of $T(A)$ if it has three neighbors in $T(A)$.
If $v_1$ and $v_2$ are the immediate descendants of $v$ and $v_3$ is the parent, then
we call $v_1$ the {\bf proper descendant} if either the edge $(v_3, v)$ points to a black face and
the path $v_3,v,v_1$ is a right turn or $(v_3,v)$ points to a white face and the edge $v_3,v,v_1$
is a left turn.  The {\bf proper branch} of $v$ is the subtree of descendants of the proper
descendant of $v$.

We say that $v \prec w$ if $v \in T_w(A)$.  The reader may easily
check that if $v$ and $w$ are adjacent, their exploration paths
agree up until the first time one of the two vertices is hit, and
thus we have either $v \prec w$ or $w \prec v$.  A similar argument
implies that $\prec$ always has a unique minimum among the vertices
in any connected subset $W$ of $V$.  In particular, for each face
$f$, there must be some minimal $v$ among the vertices of $f$.

We say a spanning tree is {\bf branch-separated} if every connected
subset of $V$ (with respect to adjacency within $G$) has a unique
minimal vertex $v$. Roughly speaking, a tree is branch-separated if
distinct branches of the tree are disconnected from one another by
the path leading from the root to the last common ancestor of the
branches. In other words, the tree can only branch at a vertex $v$
if the path from $v_0$ to $v$ passes through one of the three
neighbors of $v$ and the remaining two neighbors lie in distinct
components of the complement of that path. The above discussion
implies that exploration trees are branch separated, and the
converse is also true.

\begin{proposition} A spanning tree of $G$ is branch-separated if and only if it is the exploration
tree $T(A)$ of some $A \subset F$.  \end{proposition}

\begin{proof}
Suppose that $T$ is branch-separated and $f$ is a face of $G$.  Let
$v$ be the minimal vertex incident to $f$, and let $w$ be the
minimal vertex among the remaining vertices of $f$.  We say that $f$
is a member of $A$ if the path from $v_0$ to $w$ turns right after
hitting $v$.  The reader may check that $T= T(A)$. \end{proof}

Since the proof gives us a way to deduce $A$ from $T(A)$, we can also observe the following:

\begin{proposition} The correspondence between subsets of $F$ and branch-separated spanning trees
rooted at $v_0$ is one-to-one.  In particular, there are exactly
$2^{|F|}$ branch-separated spanning trees rooted at $v_0$.
\end{proposition}

The reader may easily check the following propositions.  We will see
analogs of these phenomena in the continuum correspondence between
loop ensembles and exploration trees.  In all of the propositions
below, unless otherwise stated, we will assume that all of the faces
outside of $G$ are colored white, so that there is a one-to-one
correspondence between subsets $A \subset F$ and collections $\{L_1,
L_2, \ldots, L_m\}$ of disjoint non-self-intersecting loops in $G$.

\begin{proposition}   The exploration tree $T(A)$ contains all but one edge of
each of the loops $L_i$. The paths in the tree traverse these loops
counterclockwise when black is on the inside and clockwise when
black is on the outside.
\end{proposition}

\begin{proposition} \label{p.onechordalpathoftreebw}
Suppose $v \not = v_0$ lies on the boundary of $G$ and that (instead
of the usual all-white coloring) the faces outside of $G$ are
colored black if they border edges in the clockwise path $P$ from
$v_0$ to $v$ in the boundary of $G$, and white if they border edges
in the counterclockwise path. Then the path $T_v(A)$ follows the
boundary between a cluster of black hexagons and a cluster of white
hexagons.
\end{proposition}

Proposition \ref{p.onechordalpathoftreebw} suggests a simple way to
construct the path $T_v(A)$ directly from the loops when $v$ is on
the boundary of $G$.  (This construction has a natural continuum
analog when $4 < \kappa \leq 8$, which is the range for which the
\CLEk/ loops intersect the boundary of $D$.)  Suppose $v \not = v_0$
lies on the boundary of $G$.  Let $M_1, \ldots, M_k$ be the loops among the
$\{L_1, \ldots, L_m\}$ that have edges in common with the directed
path $P$. Let $I_i$, for $1 \leq i \leq k$, denote the interval of
$M_k$ beginning at the first vertex of $P$ contained in $L_i$ and
ending at the last vertex in $P$ contained in $L_i$.  Let $A_i$ be
the directed arc of $M_i$ that starts and ends at the first and last
endpoints of $I_i$ and contains no edges of $P$.  Let $Q$ be the
path whose edge set is the union of the $A_i$ where $i$ ranges over
those $i$ for which $I_i$ is maximal (i.e., for which the interval
$I_i$ is not contained in any other $I_j$)
together with the edges of $P$ that are not contained in any
interval $I_i$.  Then we have the following:

\begin{proposition} \label{p.onechordalpathoftree}
The path $Q$ constructed above is equivalent to $T_v(A)$.
\end{proposition}

This in turn implies the following:

\begin{proposition} \label{p.discreteboundarytree} Let $T_{\partial G}(A)$ denote the union of $T_v(A)$
over all $v$ on the boundary of $G$.  Then $T_{\partial G}(A)$
determines---and is determined by---the set of the loops which
contain boundary edges of $G$.
\end{proposition}

%\begin{proposition} Suppose $v \not = v_0$ lies on the boundary of $G$ and that the faces outside
%of $G$ are all colored black.  Then the path from $v$ to $v_0$ in
%$T(A)$ follows the boundary of the set of all faces that are
%surrounded by loops which hit the boundary. \end{proposition}

% In particular, Loewner's theorem tells us that the path $\gamma$ can
% be reconstructed from its driving parameter.

\subsection{SLE definition} \label{s.SLEdefinition}

We will give a very brief introduction and definition of \SLE/.
There are many excellent \SLE/ surveys (most available on the arXiv) that
the reader may consult for more information on \SLE/ and the notion of scaling
limit (see, e.g., \cite{MR2079672, \LawlerConformalBook, \NienhuisKagerSurvey}).
Given a real-valued measurable function $W_t:[0,\infty) \to \R$ and
$z \in \overline \H$, consider the solution of the ODE

\begin{equation} \label{e.Loewner}
\partial_t g_t(z) = \frac{2}{g_t(z) - W_t}, \,\,\,\,\,\,\, g_0(z) =
z.
\end{equation}
When $W_t$ is continuous, and $z$ is fixed, this $g_t(z)$ is a well defined continuous function of $t$
up until the time $\inf \{t: g_t(z) = W_t \}$, at which point
it ceases to be well-defined.  More generally, if $W_t$ is merely measurable, then $g_t(z)$
is still well defined up until the first time $t$ that the set $\{g_s(z) - W_s: s < t \}$
has an accumulation point at zero.  In each case, we let $\tau_z$ be the
supremum of the $t$ for which $g_t(z)$ is
well-defined and write $K_t:= \{z \in \H: \tau_z \leq t\}$.  Then
$K_t$ is closed set, and $g_t$ is a conformal map from $\H \setminus K_t$ to $\H$.  When it
exists, we write $\gamma(t) := \lim_{z \to W_t} g_t^{-1}(z)$ (where the limit is taken
over $z \in \H$).  In the
particular case that $W_t$ is $\sqrt{\kappa} B_t$, where $B_t$ is
Brownian motion, the limit exists almost surely for all $t \geq 0$
and $\gamma: [0,\infty) \to \overline \H$ is a random continuous path,
called \SLEk/ \cite{\RohdeSchramm}.

%Given an arbitrary continuous path $\gamma$, we may denote by $f_t: \H \setminus \gamma([0,T]) \to \H$
%the conformal map which, extended to the boundary, sends $\infty$ to $\infty$ and $\gamma(T)$ to $0$.
%(This assumes that $\gamma(T)$ lies on the boundary of $\H \setminus \gamma([0,T])$.)
%Then we have the following \cite{\RohdeSchramm}:

More generally, suppose that $\gamma:[0,\infty) \to \overline \H$ is any continuous path
such that for every
$t \geq 0$ the value $\gamma(t)$ lies on the boundary of $\H \setminus \gamma([0,t])$.
For every $t\in[0,T]$, there is a unique conformal homeomorphism
$g_t:\H\setminus\gamma[0,t]$ which satisfies the so-called {\bf hydrodynamic}
normalization at infinity
$$
\lim_{z\to\infty} g_t(z)-z=0\,.
$$
The limit
$$
\capacity(\gamma[0,t]):=\lim_{z\to\infty} z(g_t(z)-z)/2
$$
is real and monotone non-decreasing in $t$.
It is called the (half plane) {\bf capacity} of $\gamma[0,t]$ from $\infty$, or
just capacity, for short.
When $\capacity(\gamma[0,t])$ is also continuous in $t$,
it is natural to reparameterize $\gamma$ so that
$\capacity(\gamma[0,t])=t$.
Loewner's theorem states that if $\gamma$ is a simple path, then the maps $g_t$ satisfy
(\ref{e.Loewner}) with $W(t)=g_t(\gamma(t))$ (where $g_t: \H \setminus \gamma[0,t] \to \H$
is extended continuously to the point $\gamma(t)$).
We can now formulate Schramm's characterization of \SLE/:

\begin{theorem} \label{conformalmarkov}
The \SLEk/ for $\kappa \geq 0$ are the only random continuous
(when parameterized by capacity) paths $\gamma:[0,\infty) \to \overline \H$
with the following so-called conformal Markov property.  Fix $T \in \R$.  Then given $\gamma$ on the
set $[0,T]$, the conditional law of $\gamma$ on the set $[T,
\infty)$ is (up to a time change) the image of the original law of $\gamma$ under the conformal
map $g_T^{-1}$ (extended continuously to $\overline \H$).  Moreover, the law of $\gamma$ is (up to a
rescaling of time) invariant under conformal automorphisms of $\overline \H$ that fix $0$ and $\infty$
(i.e., the maps $z \to az$ for $a>0$).

When $\gamma$ is an \SLEk/, the conformal Markov property continues to hold when $T$ is
replaced by an arbitrary stopping time.
\end{theorem}

It is not hard to see why this is true:
the continuity of $\gamma$ together with the conformal Markov property
implies the continuity of $W_t$.  The scale invariance implies that the law of $W_t$ is the
same as the law of $a^{-1/2}W_{a t}$ when $a>0$.
The conformal Markov property implies stationarity and independence of increments $W_{t_1} - W_{t_2}$ and
$W_{s_1} - W_{s_2}$ when $(s_1,s_2)$ and $(t_1,t_2)$ do not overlap.  These properties together imply
that $W_t$ must be a constant multiple of Brownian motion.

\subsection{$O(n)$ models and conformal invariance ansatz} \label{s.onmodel}

Let $G$ be a hexagon graph.  In the so-called {\bf $O(n)$ loop model},
one samples a collection of disjoint
non-self-intersecting loops in $G$, where each such collection has
probability proportional to $n^{N}x^{L}$ where $N$ is the number of loops, $L$ is the total
number of edges in the loops, and $n$ and $x$ are fixed positive
constants.  This can be written $\exp [\mathcal E]$ where $\mathcal
E: = N \log n + L \log x$ is a Hamiltonian on the space of loop
configurations. Equivalently, we may color all hexagons outside of
$F$ white and then sample a subset $A\subset F$ with probability
proportional to $n^{N(A)}x^{L(A)}$ where $N(A)$ is the number of
black plus the number of white clusters and $L(A)$ is the total
number of edges in the loops.  If $n=1$, then this is simply the
Ising model.  It is ferromagnetic when $x < 1$, anti-ferromagnetic
when $x
> 1$, and independent Bernoulli percolation when $x=1$.

We remark (see Figure \ref{explorationtreefigure}) that for every
loop there is exactly one edge in that loop that is not in $T(A)$;
hence $N(A) = |E_L(A) \setminus T(A)|$, where $E_L(A)$ is the set of edges
that lie in a loop, and $L(A) = |E_L|$.  Thus, we may interpret the
$O(n)$ model as a random pair $(E_L, T)=(E_L(A), T(A))$, with a Hamiltonian that
is a linear combination of $|E_L \setminus T|$ and
$|E_L|$.  (In fact, any linear combination of the four quantities
$|E_L \cup T|$, $|E_L \setminus T|$,
$|T \setminus E_L|$, and $|E_L \cap T|$ can be written this way, up to an
additive constant,
since $|T|$ and $|E|$ are both fixed.)

A natural variant of the $O(n)$ model is the following.  Fix vertices
$a$ and $b$ on the boundary of $G$ and suppose that the hexagons outside
of $F$ are colored in such a way that they are white whenever they
are incident to an edge in one arc of the boundary of
$G$ with endpoints $a$ and $b$, and black
whenever they are incident to an edge of the complementary arc. In this case, given any choice of $A$, let
$P$ be the path of edges in $F$ that lie on the boundary between the
cluster of black hexagons that includes the black boundary arc and
the cluster of white hexagons that includes the white boundary arc.
We let $L$ be the total number of edges in $E$ separating black
hexagons from white hexagons (including those edges in $P$) and let
$N$ be the total number of loops formed by these edges (not counting
the path $P$), and as before sample a subset $A\subset F$ with probability
proportional to $n^{N(A)}x^{L(A)}$.  We call this a {\bf
one-chordal-arc $O(n)$ model}.

In light of conformal invariance hypotheses from the physics
literature it is natural to
conjecture that as the mesh size gets finer, the path $P$ in the
one-chordal-arc $O(n)$ model converges in law to a random path which
satisfies a conformal Markov property---and hence, by Theorem \ref{conformalmarkov},
is \SLEk/ for some $\kappa$.
Neinhuis and Kager \cite{\NienhuisKagerSurvey}, following work by
Duplantier, Neinhuis, and others (which uses the so-called
``Coulomb gas method''),
have gone further and conjectured precise values for $\kappa$: namely, they conjecture that the
scaling limit of $P$ is
\begin{enumerate}
\item \SLEk/, where $n = -2\cos (4 \pi/\kappa)$ and $4 \leq \kappa \leq 8$, if $0 < n \leq 2$
and $x > x_c$,
\item \SLEk/, where $n = -2\cos (4 \pi/\kappa)$ and $8/3 \leq \kappa \leq 4$, if
$0 < n \leq 2$ and $x = x_c$, \item a straight line (or a shortest
length path from $a$ to $b$, if the domain is not convex) if either
$x < x_c$ or $n > 2$,
\end{enumerate}
where $x_c = [2 + (2-n)^{1/2}]^{-1/2}$.  Note that $x_c$ increases
monotonically as $n$ increases from $0$ to $2$, and for each $n \in
(0,2)$, the equation $n = -2\cos (4 \pi/\kappa)$ has two solutions,
one in $(8/3,4)$ and one in $(4,8)$.  A precise version of this
conjecture has been proved in the case $x_c = 1/2$ and $n=1$ (which
corresponds to critical Bernoulli site percolation)
\cite{\CamiaNewman,\SmirnovPerc}.

\subsection{Excursions and renewal times} \label{s.excursion}

We now describe one way to construct loops from trees---a continuum
version of which will be used in Section \ref{CLEconstruction} to
define the \CLEk/.
Fix a vertex $v$ of $G$.  For each $k$, let $K_k$ be the set of
hexagons whose colors are determined by the first $k$
vertices $v_1, v_2, \ldots v_k$ of the exploration path $T_v(A)$.
Let $G_k$ be the connected component of the set of faces of $G$ in
the complement of $K_k$ that contains a face incident to $v$.  The
reader may observe that for each $k$, the conditional law of the
coloring of the faces within $G_k$, given the colors of the faces
determined by $v_1, \ldots, v_k$, is given by either an $O(n)$ model
(i.e., the hexagons on the boundary of $G_k$ are all one color) or a
one-chordal-arc $O(n)$ model (i.e., boundary conditions are given by
one white and one black arc).  Let $0=k_0, k_1, k_2, \ldots$ be the
values of $k$ for which the boundary of $G_k$ is monochromatic. We
refer to these $k_i$ as {\bf renewal times} and to the paths between
these times consisting of edges that separate black and white faces as {\bf excursions}.  The reader
may observe (see Figure \ref{explorationtreefigure}) that each
excursion traces part of a loop.  If we extend the excursion to a longer
directed path in $T(A)$ in such a way that the path turns in the proper direction
at each branch point of $T(V)$, then this extended path will trace out
the remainder of the loop (except for one edge).

\subsection{Two exploration path variants}

We now describe two variants of the exploration path.  In later
sections we will use continuum analogs of these variants to define
\CLEk/ when $\kappa \leq 4$.  In the exploration tree described
above, the orientations of the outermost loops are all clockwise,
and the sequence of nested loops surrounding a single face strictly
alternates between clockwise and counterclockwise.

\begin{flushleft} {\bf Variant 1:} Fix a parameter $\beta
\in [-1,1]$ and independently orient each outermost loop
clockwise with probability $\frac{1-\beta}{2}$ and counterclockwise
with probability $\frac{1+\beta}{2}$.  Inductively, we define
orientations for the remaining loops, orienting the loop the same
way as the smallest (in terms of enclosed area) loop that surrounds it with probability
$\frac{1-\beta}{2}$ and the opposite way with probability
$\frac{1+\beta}{2}$. We also independently assign an orientation to each isolated
vertex (i.e., a vertex which does not lie in any loop), orienting
the vertex the same way as the smallest loop that surrounds it
(or clockwise if the vertex is not surrounded by a loop) with
probability $\frac{1-\beta}{2}$.  The case $\beta = 0$ is
particularly natural, since in this case all loops and vertices are
oriented independently with fair coins.
\end{flushleft}
Given an oriented collection of loops and vertices, we can define an
{\bf oriented exploration tree} as follows: Fix a vertex $v_0$ on the outer
boundary of $G$ which is incident to only two edges of $G$, and fix
a directed edge $e_0$ outside of $G$, beginning at some vertex
$v_{-1}$ and pointing towards $v_0$. For each vertex $v$ of $G$, the
{\bf oriented exploration path} $T_v(A)$ is a directed,
non-self-intersecting path $v_0, v_1, v_2 \ldots$ that ends at $v$.
Each $v_k$ is chosen in such a way that the sequence $v_{k-2},
v_{k-1}, v_k$ describes a left turn if $v_{k-1}$ is a
counterclockwise isolated vertex, a right turn if $v_{k-1}$ is a
clockwise isolated vertex, and a turn that causes $(v_{k-1}, v_k)$ to
be an edge of a loop oriented in the direction of the loop if $v$ is
not isolated --- {\em unless} this choice of $v_k$ would fail to lie in
the same connected component of $V \backslash \{v_0, v_1, \ldots,
v_{k-1} \}$ as $v$, in which case the path turns the other
direction.  The reader may check that the union of the $T_v(A)$ is a
spanning tree.  We call this the oriented exploration tree.
Note that our original definition of exploration tree corresponds to
the case $\beta = 1$, while the mirror image (i.e., the model with the roles of
black and white reversed) corresponds to $\beta =
-1$.

\begin{flushleft} {\bf Variant 2:} To motivate this variant, we note that we will eventually want
to conjecture the existence of a scaling limit for the individual
exploration tree paths $T_v(A)$ which has a continuous Loewner
evolution. However, when $\kappa \leq 4$ and $\beta = 1$, we cannot
expect the scaling limit of $T_v(A)$ to have a continuous Loewner
evolution when $v$ is on the boundary of $G$, since we expect there
to be no macroscopic loops that hit the boundary; in this case, in light of
Proposition \ref{p.onechordalpathoftreebw}, we would expect $T_v$ to tend
to a path that traces the left boundary of the domain from $v_0$ to
$v$.
\end{flushleft}
In the $O(n)$ model setting, recall from the previous section that
given the values $v_0, v_1, \ldots, v_k$ up to a renewal time $k$,
the law of the remainder of the path $T_v(A)$ is that of an
exploration path from $v_k$ to $v$ in a new hexagon graph $G_k$ with
a new starting point $v_k$.  Now suppose that we generate $T_v(A)$
as above except that at each renewal time $k$, we ``shift'' this new
starting point by some constant-order number of edges to the right
along the boundary of $G_k$ (adding these edges to $T_v(A)$) before
continuing to grow $T_v(A)$ according to the usual rules. Then we
get a variant of $T_v(A)$ which we might expect (if the shift sizes
are chosen properly) to have a meaningful scaling limit. We will
make these notions more precise in the continuum setting, where the
shifts will be replaced with a local time L\'evy compensation used to make
\SLEr \kappa-6/ well-defined when $8/3 < \kappa \leq 4$.

\section{Constructing \SLEkr/ from Bessel and stable processes}
\label{s.besselsection}

\subsection{Bessel processes}
In this section we define Bessel processes and state some standard
facts that will be useful in defining conformal loop ensembles. See
Chapter XI of \cite{\RevuzYor}, including Exercises 1.25 and 1.26.
The {\bf square Bessel process} of dimension $\delta>0$, written
\BESQd/, is the unique strong solution $Z_t$ to the SDE:

\begin{equation} \label{e.besqsde} Z_t = Z_0 + 2 \int_0^t \sqrt{Z_s} dB_s + \delta
t,
\end{equation}
where $Z_0$ is a fixed initial value and $B_t$ is standard Brownian
motion.  (Recall that a strong solution to (\ref{e.besqsde}), as
defined e.g. in  \cite{\RevuzYor}, is a coupling of $Z_t$ and $B_t$
in which (\ref{e.besqsde}) almost surely holds for all $t$ and $Z_t$
is adapted to the filtration generated by $B_t$---i.e., the law of
$Z_t$ is non-anticipative.)

The {\bf Bessel process} of dimension $\delta$, written \BESd/, is
the process $X_t = \sqrt{Z_t}$, where $Z_t$ is a \BESQd/. We
sometimes write \BESdx/ for the \BESd/ process started at $X_0 =
x\geq 0$.

\begin{proposition}\label{p.zerolebesgue} If $\delta \geq 2$, and $X_t$ is a \BESd/,
then almost surely $X_t>0$ for all $t > 0$. If $0 < \delta < 2$,
then $X_t$ almost surely assumes the value zero on a non-empty
random set with zero Lebesgue measure.
\end{proposition}

\begin{proposition}\label{p.brownianscaling}
For all $\delta > 0$, the Bessel processes $X_t$ are invariant under
Brownian scaling.  That is, given a constant $c>0$, the process
$c^{-1/2} X_{c t}$ has the same law as $X_t$.
\end{proposition}

\begin{proposition}\label{p.semimartingale} If $\delta > 1$ and $X_t$ is a \BESd/ process, then
$X_t$ is a semimartingale and a strong solution to the SDE
\begin{equation} \label{BesselSDE}
X_t = X_0 + B_t + \frac{\delta-1}{2} \int_0^t X_s^{-1} ds.
\end{equation}
If $\delta \leq 1$ and $X_0>0$, then (\ref{BesselSDE}) still holds
up until the first $t$ for which $X_t = 0$, but for all larger $t$,
the integral $\int_0^t X_t^{-1} ds$ is infinite (so that
(\ref{BesselSDE}) cannot hold).  If $\delta > 0$ and $X_t$ is any
continuous process adapted to the filtration generated by $B_t$
which is {\bf instantaneous reflecting} at zero (i.e., the Lebesgue
measure of $\{t: X_t = 0 \}$ is almost surely zero) and almost surely satisfies
$$\frac{\partial}{\partial t}(X_t - B_t) =\frac{\delta-1}{2}
X_t^{-1}$$ whenever $X_t \not = 0$, then the law of $|X_t|$ is that of
a \BESd/ process.
\end{proposition}

\begin{proposition}
When $\delta = 1$ and $X_t$ is a \BESd/, the process $X_t$ has the
law of the absolute value of a standard Brownian motion.  It also
satisfies the equation

\begin{equation} \label{BesselSDE1}
X_t = B_t + \frac{1}{2} l^0_t,\end{equation} where $l^0$ is the
local time of $X_t$ at $0$, where the local time $l^x_t$ is defined
to be the almost surely continuous function of $x$ and $t$ for which
$$\int_0^t f(X_s)ds = \int_0^\infty f(x)l^x_t dt,$$ for all $t>0$
and measurable functions $f$.
\end{proposition}

\begin{proposition} \label{p.occupationdensity}
When $\delta \in (0,1)$, there almost surely exists a family of occupation densities
$l^x_t$ for a Bessel process $X_t$ that are continuous in $x$ and $t$
such that for each $t>0$ and measurable $f$, we have

$$\int_0^t f(X_s)ds = \int_0^\infty f(x) l^x_t x^{\delta-1}dx.$$
For $t>0$ and $\delta \in (0,1)$, the integral $\int_0^t X_s^{-1}ds$
is almost surely infinite.  However, the so called {\bf principal
value}
$$\PV \int_0^t X^{-1}_s ds := \int_0^\infty x^{\delta - 2}(l^x_t -
l^0_t)ds$$ is almost surely finite for all $t$, and satisfies

\begin{equation} \label{BesselSDEsmalldelta}
X_t = X_0 + B_t + \PV \frac{\delta-1}{2} \int_0^t X_t^{-1} ds.
\end{equation}

\end{proposition}

\medskip

In light of (\ref{BesselSDEsmalldelta}) we may take $\PV \int_0^t
X_t^{-1} ds = \frac{2}{\delta-1} (X_t -X_0- B_t)$ as an alternate
definition for the principal value when $\delta \in (0,1)$ (when the
coupling of $X_t$ and $B_t$ is given).  It is clear from
(\ref{BesselSDEsmalldelta}) that this integral is also a process
that satisfies Brownian scaling.

\subsection{L\'{e}vy skew stable processes}
\label{s.levy}

We now review some basic facts about L\'evy skew stable
distributions and their connection to Bessel processes and skew
Bessel processes. They are not hard to derive directly, but the
reader may see, e.g.,
\cite{\PitmanPart,\RevuzYor,\SkewBesselSurvey,\Watanabe} and the
references therein for more details and many additional results. The
L\'{e}vy skew stable probability distribution is the
Fourier transform of its characteristic function $\phi$, defined as
follows. Fix parameters $c> 0$, $\beta \in [-1,1]$, $\mu \in \R$ and $\alpha \in
(0,2]$. Then

\begin{equation}\label{e.skewstable}\phi(\lambda) = \exp \left[i\lambda\mu -|c\lambda|^\alpha(1-i\beta \sign(\lambda) \Phi)
\right],\end{equation} where $\Phi = \tan(\pi \alpha/2)$ if $\alpha \not =
1$, and $\Phi = -(2/\pi) \log|\lambda|$ for $\alpha = 1$.  When $\alpha \in (0,2)$, the
corresponding L\'{e}vy measure is
$$\Lambda(d\eta) = \frac{c \alpha}{\Gamma(1 - \alpha)}
\eta^{-\alpha-1}d\eta \,\,\,\,\,\,\,\, (t > 0),$$
if $\beta = 1$ and more generally the measure $\Lambda_\beta$ defined by
$$\Lambda_\beta(A) := \frac{1+\beta}{2} \Lambda(A) + \frac{1-\beta}{2} \Lambda(-A),$$
for measurable subsets $A$ of $\R$.

For each choice of parameters $\alpha, \mu, \beta$ as above and constant $b>0$,
there is a corresponding stable process $S_t$ with independent,
stationary increments, such that for each fixed $t$, the law of
$S_t$ is given by the L\'evy skew stable distribution with parameter
$c$ where $c^\alpha=b t$.  We denote this process by $S(\alpha, \beta, \mu, b)$.
It is supported on the positive
reals if and only if $\mu \geq 0$, $\beta = 1$, and $\alpha < 1$.  The jump discontinuities
in $S_t$ have sizes that are distributed as a Poissonian point process sampled from the
corresponding L\'evy measure.  When $b=1$, this means in particular that the expected number
of jump discontinuities in $S_t$ whose size lies in a set $A$ that occur between time $s_1>0$ and
time $s_2>s_1$ is given by $(s_2-s_1)\Lambda_\beta(A)$.  When $\beta=1$ the jumps are almost
surely all positive; when $\beta=-1$ they are almost surely all negative.

The following is
easy to derive from the Markov and scaling properties of Bessel
processes:

\begin{proposition} \label{p.localtimeskewstable}
Let $l = l(t) = l^0_t$ be the zero local time of a Bessel process
$X_t$ with parameter $\delta \in (0,1) \cup (1,2)$. Then $t(l)$ and
$\PV \int_0^{t(l)} X_s ds$ are both L\'evy skew stable processes
indexed by $l$ with parameters $\beta =1$, $\mu = 0$, $\alpha = 1-\delta/2$ and
$\alpha = 2-\delta$ respectively, and some positive $b$.  When
$\delta \in (0,2)$, the zero set of a Bessel of dimension $\delta$
is the range of a non-decreasing stable process (a.k.a. stable
subordinator) with parameter $\alpha = 1-\delta/2$.
\end{proposition}

We say that an $S(\alpha, \beta, \mu, b)$ process is {\bf strictly stable} if altering
$b$ (i.e., rescaling time by a constant factor) has the same effect on the law
of the process as multiplying the process by a deterministic constant.
The following is not hard to derive from (\ref{e.skewstable}):

\begin{proposition}
Fix $\alpha \in (0, 2)$, $\beta \in [-1,1]$, and $\mu \in \R$.  Then
the L\'evy skew stable process $S(\alpha, \beta, \mu, b)$ is
strictly stable if and only if one of the following holds:
\begin{enumerate}
\item $\alpha \not = 1$, $\mu = 0$.
\item $\alpha = 1$, $\beta = 0$.
\end{enumerate}
\end{proposition}

We remark that it is easy to see from (\ref{e.skewstable}) that one
can obtain the characteristic function corresponding to $\alpha = 1$, $\beta = 0$ and $\mu \not = 0$
as a uniform limit of the characteristic functions corresponding to $\alpha \not = 0$, $\mu = 0$ and
$\beta \not = 0$ if one takes $\beta \to 0$ at an appropriate rate as $\alpha \to 1$.

The {\bf skew Bessel process} with dimension $\delta$ and skew
parameter $\beta$ is a continuous process $X_t$ for which the law of
$|X_t|$ is that of the Bessel process of the corresponding
dimension, but for each excursion of $|X_t|$ (i.e., each connected
component of $\{t:|X_t|>0\}$) we toss an independent coin and make
$X_t$ positive on that excursion with probability
$\frac{1+\beta}{2}$ and negative with probability
$\frac{1-\beta}{2}$.  When $\beta=0$, the resulting process is
called the symmetric Bessel process.  Now we define a process $Y_t$
that will play the same role as the principal value $\PV
\int_0^{t(l)} X_s ds$ when $\beta \not = 1$.

\begin{proposition} \label{p.localtimeskewstableotherbeta}
Let $l = l(t) = l^0_t$ be the zero local time of a skew Bessel
process $X_t$ with parameters $\beta \in [-1,1]$ and $\delta \in
(0,2)$.  If either $\delta \not = 1$ or $\delta=1$ and $\beta=0$,
then $X_t$ can be coupled with a continuous process $Y_t$ such that
\begin{enumerate}
\item The pair $(X_t, Y_t)$ is adapted to the filtration generated
by the Brownian motion $B_t$.
\item $Y_t' = X_t^{-1}$ on the set $\{t: X_t \not = 0 \}$, almost surely.
\item The law of $(X_t,Y_t)$ is invariant under Brownian scaling.
\item If $T$ is a stopping time of $(X_t, Y_t)$ for which $X_T = 0$
almost surely, then $T$ is a renewal time in that conditioned on
$T$, the law of $(X_{T+t}, Y_{T+t} - Y_T)$ (for $t \geq 0$) is the
same as the original law of $(X_t, Y_t)$.
\end{enumerate}
In any such coupling, $Y_{t(l)}$ is an $S(2 - \delta, \beta, \mu,
b)$ process indexed by the local time parameter $l$, where $b$ is as
given in Proposition \ref{p.localtimeskewstable} (for the case
$\beta=1$) and $\mu = 0$ unless $\delta = 1$. The law of $(X_t,Y_t)$
is uniquely determined by the properties above (together with the
parameter $\mu$, in the case $\delta = 1$).
\end{proposition}

It is not difficult to derive Proposition
\ref{p.localtimeskewstableotherbeta} from Proposition
\ref{p.localtimeskewstable} when $\delta \not = 1$.  We sketch
the construction of $Y_t$ as follows.  To construct $Y_t$, it is
sufficient to determine the values of $Y_t$ on the set $\{t: X_t = 0
\}$, since the other values for $Y_t$ may then be obtained by integrating
$X_t^{-1}$ on each excursion.  Thus it is enough to determine the process $Y_{t(l)}$
as a function of $l$. When
$\beta = 1$, the jump discontinuities in $Y_{t(l)}$ (each of which
corresponds to the integral of $X_t^{-1}$ over an interval on which
$X_t \not = 0$) are distributed according to the Poisson process
derived from the L\'evy measure $\Lambda$ described above; swapping
the sign of a $\frac{1-\beta}{2}$ fraction of these jumps
corresponds to replacing $\Lambda$ with $\Lambda_\beta$ where
$$\Lambda_\beta(A) = \frac{1+\beta}{2}\Lambda(A) + \frac{1-\beta}{2}
\Lambda(-A),$$ which in turn corresponds to changing the skew
parameter of the process $Y_{t(l)}$ from $1$ to $\beta$.  When
$\delta = 1$, $\mu = 0$, and $\beta = 0$, the process $(X_t, Y_t)$
may be obtained as a limit as $\delta \to 1$ of the $(X_t, Y_t)$
couplings with $\beta=0$ and $\mu = 0$.  Adding a non-zero value for
$\mu$ amounts to replacing $Y_t$ with $Y_t + \mu l^0_t$.

\subsection{Chordal \SLEkr/} \label{s.chordalsle}

Fix a constant $\rho \in R$.  Write $\delta = 1 +
\frac{2(\rho+2)}{\kappa}$.  Suppose $\delta > 0$ and $\delta \not =
1$ (i.e., $\rho \not = -2$).  Let $X_t$ be a \BESdx/, and let $O_t$
and $W_t$ be given by

\begin{eqnarray*}
O_t &=& -2 \kappa^{-1/2} \PV \int_0^t X_s^{-1}ds\\
W_t &=&  O_t + \sqrt \kappa X_t\\
\end{eqnarray*}
where initial values $W_t = \sqrt \kappa x$ and $O_t < W_t$ are
given.  We then define \SLEkr/ to be the growing family of closed
sets $K_t$ determined by the Loewner evolution with the driving
parameter $W_t$ as given above.  (See the definition of \SLE/,
Section \ref{s.SLEdefinition}.)

More generally, let $\beta \in [-1,1]$, $\kappa>0$, $\mu \in \R$,
and $\rho \in \R$ be given and define $\delta = 1 +
\frac{2(\rho+2)}{\kappa}$ as above. If either $\delta \in (0,1) \cup
(1,2)$ and $\mu=0$ or $\delta = 1$ and $\beta = 0$, then we define
{\bf skew \SLEkr/ with parameters $\beta$ and $\mu$}, which we
denote \SLEbkr/, the same way as we defined \SLEkr/ above except
that we replace $\PV \int_0^t X_s^{-1}$ with the process $Y_t$ of
Proposition \ref{p.localtimeskewstableotherbeta}.  In other words,
we begin with the pair $(X_t,Y_t)$ from Proposition
\ref{p.localtimeskewstableotherbeta} (with some initial values $X_0$
and $Y_0$ fixed) and define $O_t = -2 \kappa^{-1/2} Y_t$ and $W_t =
O_t + \sqrt \kappa X_t$.  Note that \SLEkr/ is equivalent to
\SLEbkr/ with $\beta = 1$, $\mu = 0$.  We will assume $O_0=W_0=0$
when we don't specify otherwise.

%The following is immediate from the fact that the pair $(X_t, Y_t)$
%satisfies Brownian scaling (Proposition
%\ref{p.localtimeskewstableotherbeta}).
%\begin{proposition}
%The processes \SLEkr/ and \SLEbkr/ defined above are scale
%invariant---that is, when $c > 0$ is fixed the process $c K_t$ has
%the same law as $K_t$, up to a constant time change.
%\end{proposition}

Given a domain $D$ with marked boundary points $a$ and $b$, a chordal \SLEk/
from $a$ to $b$ in $D$ is the image of chordal \SLE/ from $0$ to
$\infty$ in $\H$ (as defined above) under any conformal map taking
$a$ to $0$ and $b$ to $\infty$.  We will be particularly interested
in chordal \SLEkk \kappa; \kappa-6/ because of the following
\cite{\SchrammWilson}:

\begin{proposition} \label{p.chordalinvariance} In the chordal \SLE/ context, fix $W_0=0$ and $O_0 = a$ for
some non-zero $a \in \R$.  Then an \SLEr \kappa - 6/ (from $0$ to
$\infty$) in $\H$, with these initial values --- stopped at the
first time $W_t = O_t$
--- has the same law (up to time change) as a chordal \SLEk/ in $\H$
from $0$ to $a$ (which we may view as a random path $\gamma$)
--- stopped the first time $t$ for which $a$ and $\infty$ fail to lie on the boundary of the same
connected component of $\H \setminus \gamma([0,t])$.
\end{proposition}

The reader may view the following as a reason to expect \SLEr
\kappa-6/ to be the scaling limit of the $T_v(A)$ described in
Section \ref{s.discretesection} when $A$ is the coloring
corresponding to an $O(n)$ model and $v$ is a boundary vertex.

\begin{proposition} \label{p.Ktinvariance}
Suppose $K_t$ is a random Loewner evolution in $\R$, driven by some
continuous $W_t$, and write $O_t = g_t (\inf \{K_t \cap \R\})$.
Suppose that $K_t$ satisfies the following:
\begin{enumerate}
\item {\bf Scale Invariance:} for any positive constant $c$,
the law of the process $t \to K_t$ is the same as that of the
process $t \to cK_t$ up to time parameterization.
\item {\bf Renewal property:} given $W_t$ up to any stopping time $T$ for
which $O_T = W_T$ almost surely, the conditional law of the process $W_{T+t} -
W_T$, for $t \geq 0$, is the same as the original law of $W_t$.
\item {\bf Conformal Markov property:} given $W_t$ up to any fixed time $T$,
the conditional law of the growth process $\overline{g_T K_{T+t}}$ --- up to
time $\inf \{t : t \geq 0, O_{T+t} = W_{T+t} \}$ --- is the same as that of an
ordinary chordal \SLEk/ in $\H$ from $W_T$ to $O_T$ up to that time.  (Here $g_T K_{T+t}$ is
a subset of $\H$ but the closure is taken in $\overline \H$.)
\end{enumerate}
Then $K_t$ is an \SLEr \kappa-6/ for some $\kappa > 4$.  Conversely, the
three properties above hold more generally when $W_t$ and
$O_t$ are as in the definition of the \SLEbr \kappa-6/ process with
$\kappa > 8/3$, $\mu \in \R$, and $\beta \in [-1,1]$ (provided
$\beta=0$ if $\kappa=4$ and $\mu = 0$ if $\kappa \not = 4$),
although in this generality it is no longer the case that $O_t =
g_t(\inf \{K_t \cap \R\})$.  The conformal Markov property holds when
$T$ is replaced with an arbitrary stopping time $T$ for which $O_T \not = W_T$ almost surely.
\end{proposition}

\begin{proof}
First we claim that $X_t = \kappa^{-1/2} (W_t-O_t)$ is a Bessel
process with $\delta = 1 + \frac{2(\rho+2)}{\kappa}$, where $\rho =
\kappa-6$.  It follows from Proposition \ref{p.chordalinvariance}
that when $X_t > 0$, it evolves according to the SDE for this
process.  Since $X_t$ is almost surely positive, it is sufficient by
Proposition \ref{p.semimartingale} to show that the Lebesgue measure
of the set of times $t$ for which $W_t= g_t (\inf \{K_t \cap \R\})$
is almost surely zero. In fact, we claim that this is true for any
continuous Loewner evolution $W_t$.

To see this, fix $\epsilon > 0$ and let $A_t$ be the largest integer multiple of
$\epsilon$ less than $\inf \{K_t \cap \R\}$.  The process
$\tilde O_t = g_t(A_t)$ evolves differentiably according to (\ref{e.Loewner}) except at
discrete times when it jumps by discrete amounts to the left, and
the monotonicity of (\ref{e.Loewner}) implies that $|\tilde O_t - O_t|
\leq \epsilon$ for all $t$. Now, $(W_t - O_t) < \epsilon$ implies
$W_t - \tilde O_t < 2 \epsilon$. Let $T$ be the first time for which $O_t
\leq C$, for some constant $C < 0$. Since $\tilde O_0 = -\epsilon$ and
$\tilde O_T \geq C - \epsilon$, it follows from (\ref{e.Loewner}) that the
Lebesgue measure of $\{t: 0 \leq t \leq T, |W_t - O_t| < \epsilon\}$
is less than or equal to $\epsilon C$. Since this holds for any
$\epsilon$, it in particular implies that $\{t: 0 \leq t \leq T, W_t
= O_t\}$ has Lebesgue measure zero. Since this holds for any $C$, it
proves the claim.

Second, the reader may easily check that $X_t = \kappa^{-1/2}
(W_t-O_t)$ and $Y_t = -2 \sqrt \kappa O_t$ satisfy the hypotheses of
Proposition \ref{p.localtimeskewstableotherbeta}; the fact that
$Y_t' = X_t^{-1}$ on the set $\{t: X_t \not = 0 \}$, almost surely,
follows from Proposition \ref{p.chordalinvariance} and the conformal
Markov property, while the other properties are consequences of the
scale invariance and renewal assumptions.  It follows that the
$(O_t, W_t)$ is the pair arising in some \SLEbkr/.  The fact that
$O_t$ is almost surely non-decreasing implies that $\beta = 1$ and
$\delta > 1$; hence $\kappa > 4$.

The concluding two sentences of Proposition \ref{p.Ktinvariance} are
immediate from Propositions \ref{p.chordalinvariance} and
\ref{p.localtimeskewstableotherbeta}.
\end{proof}

An ordinary \SLE/
starting from a boundary point $a$ and ending at a boundary point
$b$ on a planar domain is believed to have the same law (up to a
time change) as an \SLE/ starting at $b$ and ending at $a$
\cite{\RohdeSchramm}.  (This must be the case for $8/3 < \kappa < 8$ if \SLEk/ is
the scaling limit of the one-boundary-arc $O(n)$ models discussed in
Section \ref{s.onmodel}.) However, this invariance does not readily
follow from the definition of \SLE/.  The following conjecture will turn
out to be relevant to the study of conformal loop ensembles:

\begin{conjecture} \label{c.reversible}
Fix $4 < \kappa < 8$ (and $\mu = 0$, $\beta = 1$).  Then the processes \SLEk/ and chordal \SLEkk
\kappa;\kappa - 6/ are both almost surely continuous paths.  The laws
of these paths---up to direction of parameterization---are invariant
under anticonformal automorphisms of $D$ that swap the endpoints
of the paths.
\end{conjecture}

The fact that \SLEk/ is continuous appears in \cite{\RohdeSchramm}, but the proof
has never been extended to \SLEkr/ processes.

\subsection{Radial \SLEkr/} \label{s.radialSLEkr}
We will now introduce radial \SLEkr/.  Let $\D \subset \C$ be the
unit disc. Radial \SLEk/ is a random path from a boundary point of
$\D$ to the center of $\D$ defined the same way as chordal \SLEk/ except
that instead of (\ref{e.Loewner}) we use the ODE$$\partial_t g_t(z) = \Psi(W_t, g_t(z)),$$
where $W_t$ is a point on the unit circle and (following notation
from \cite{\SchrammWilson})
$$\Psi(w,z) := -z \frac{z+w}{z-w}.$$ In this case \SLEk/ is defined by taking
$W_t = e^{i \sqrt \kappa B_t}$ where $B_t$ is a standard Brownian
motion. Equivalently, $W_t$ is the solution to the SDE
$$dW_t= (- \kappa/2) W_t dt + i \sqrt \kappa W_t dB_t.$$
Given any measurable driving function $W_t:[0,\infty) \to \partial \D$,
we let $\tau_z$ be the supremum of the $t$ for which $g_t(z)$ is
well-defined and write $K_t:= \{z \in \D: \tau_z \leq t\}$.  Then
$g_t$ is a conformal map from $\D \setminus K_t$ to $\D$.  When it
exists, we write $\gamma(t) := \lim_{z \to W_t} g_t^{-1}(z)$, as in
the chordal case.  (Recall that in the setting of radial \SLEk/, $\gamma$ exists
and is continuous almost surely for all $\kappa \geq 0$ \cite{\RohdeSchramm}.)
Note that $t$ represents not the half-plane capacity
of $K_t$ (as in the chordal case) but $-1$ times the log of the conformal
radius of $\D \setminus K_t$ viewed from zero (i.e., $t=\log|g_t'(0)|$).

If we take $O_t$ to be another point on the unit circle, we can
define radial \SLEkr/ --- at least up until the first time $O_t$ and
$W_t$ collide --- by taking

$$dO_t = \Psi(W_t, O_t)dt$$
and
$$dW_t = (- \kappa/2) W_t dt + i \sqrt \kappa W_t dB_t +
\frac{\rho}{2}\tilde\Psi(O_t, W_t)dt,$$ where $$\tilde \Psi(z,w) :=
\frac{\Psi(z,w) + \Psi(\overline z^{-1}, w)}{2}.$$  Given initial
values $O_0 \not = W_0$ the solution to this SDE exists uniquely up
until the first time $W_t = O_t$; see \cite{\SchrammWilson}, which
also proves the following analog of Proposition
\ref{p.chordalinvariance}:

\begin{proposition} \label{p.slecoordinatechange1}
Fix $W_0 = 1$ and $O_0 = a$ to be distinct points on the unit circle
$\partial \D$.  Then a radial \SLEr \kappa - 6/ in $\D$ (from $1$ to
$0$) --- stopped the first time $W_t = O_t$
--- has the same law (up to time change) as a chordal \SLEk/
path $\gamma$ in $\D$ from $1$ to $a$
--- stopped the first time $t$ for which $a$ fails to lie
on the component of $\H \setminus \gamma([0,t])$ containing $0$.
\end{proposition}

We can extend the definition of radial \SLEr \kappa-6/ beyond times
for which $W_t = O_t$ by mapping the corresponding chordal \SLEr
\kappa-6/ into the unit disc.  If $O_t$ and $W_t$ are continuous processes
on $\partial D$, then let $\hat O_t$ denote the lifting of $\arg (W_t - O_t)$
to a continuous function on $\R$.  We will prove the following in Section
\ref{s.approx}:

\begin{proposition} \label{p.radialchordalcoordchange}
Fix $\kappa \in (8/3, 8)$, $\beta \in [-1,1]$, and $\mu \in \R$ (such that
$\beta = 0$ if $\kappa = 4$ and $\mu = 0$ if $\kappa \not = 4$).  Then there exists a
unique continuous Markovian diffusion on pairs $(W_t, \hat O_t)$ with
the following property: when the initial values $(W_0, \hat O_0)$ are such
that $\hat O_0 = 2 k \pi$ for some
integer $k$ and $\psi$ is any conformal map (if $k$ is even) or an anti-conformal map
(if $k$ is odd) from $\D$ to $\H$ for which $\psi(W_0)=0$, the image of
the corresponding radial Loewner evolution $K_t$ under $\psi$---up until
the first time $t$ that $\psi^{-1}(\infty) \in K_t$---is given by
chordal \SLEbr \kappa-6/ with initial values $W_0=O_0=0$ (up until that time).
\end{proposition}

We then define {\bf radial \SLEbr \kappa-6/} (in $\D$ with target $0$) to be the radial Loewner
evolution driven by the $W_t$ from Proposition \ref{p.radialchordalcoordchange}.  When
not specified otherwise, we take initial values to be $W_0=1$ and $\hat O_0 = 0$ (so that
$O_t = 1$).  We can then define radial \SLEbr \kappa-6/ in any other domain $D$ with a fixed
target $z \in D$ to be the image of the radial \SLEbr \kappa-6/ defined above under a conformal map
that sends $0$ to $z$ and $W_0$ and $O_0$ to the appropriate initial values on the boundary
of $D$.  The following is an immediate consequence of Proposition \ref{p.radialchordalcoordchange}
(and the fact that the choice of $\phi$ in the statement was arbitrary).

%The Loewner evolution generated by the $W_t$ above
%is called {\bf radial \SLEr \kappa-6/}.  We can use the
%r which $a_1$ fails to lie on the boundary of the component $\D_t$
%of $\D \setminus K_t$ containing $0$.  Then we pick a new point
%$a_2$ on $\partial \D_{t_1}$ and continue with a new \SLEr
%\kappa-6/ process, with the initial values of $O_{t_1}$ and
%$W_{t_1}$ as given (if $a_2$ lies on the clockwise arc from
%$O_{t_1}$ to $W_{t_1}$) or the mirror image of an \SLEr \kappa-6/
%process (if $a_2$ lies on the counterclockwise arc from $O_{t_1}$ to
%$W_{t_1}$).
%It is not hard to show that if the $a_i$ are chosen appropriately,
%the resulting growth process $K_t$ will eventually reach arbitrarily
%closely to $0$.  (We have not yet proved that this definition is
%independent of the way the $a_i$ are chosen.) We can define radial
%\SLEbr \kappa-6/ similarly for all the values of $\kappa$, $\mu$,
%and $\beta$ for which chordal \SLEbr \kappa-6/ was defined. We will
%give a formal and complete definition in Section
%\ref{s.approxslekr}, where we will also prove the following:

\begin{proposition} \label{p.invariance}
For any $\kappa > 8/3$, $\mu \in \R$, and $\beta \in [-1,1]$ (where $\beta = 0$ if
$\kappa=4$, $\mu=0$ if $\kappa \not = 4$) the law of radial \SLEbr \kappa-6/ is target invariant.
That is, if we fix initial values $W_0=O_0=1$ and $\hat O_0=0$ and we fix
distinct points $a,b \in \D$, then the law of \SLEbr \kappa-6/
targeted at $a$ and the law of \SLEbr \kappa-6/ targeted at
$b$
--- both defined up to supremum of the set of times $t$ for which $a \not \in K_t$
and $b \not \in K_t$ --- are the same (up to a time change).
\end{proposition}

\section{Conformal loop ensembles}
\label{CLEsection}
\subsection{Defining loops}

Before defining \CLE/, we need to define a suitable space of loops.
(See also \cite{\AizemanBurchard,\CamiaNewman} where similar spaces of
loops and $\sigma$-algebras on loop ensembles are introduced.)
A simple loop is a subset of $\C$ which is homeomorphic to the
circle $S^1$. Equivalently, we identify a simple loop with a
homeomorphism from $S^1$ to a subset of $\C$, modulo monotone reparameterization.

Let $\mathcal C$ be the set of homeomorphisms from $S^1$ to
subsets of $\C$ and let $\overline{\mathcal C}$ denote the closure of
$\mathcal C$ with respect to the $L^\infty$ metric.  A {\bf
quasisimple loop} is an element of $\overline{\mathcal C}$, modulo monotone reparameterization.
A quasisimple loop may
intersect itself, but it cannot ``cross'' itself transversely.
Clearly, the winding number of a quasisimple loop $L$ around a point
$z \in \C \setminus \eta(S^1)$ is equal to zero or $\pm 1$ (since
this is true of simple loops and remains true under uniform limits
provided that the limit fails to intersect $z$).
In the latter case, we say that $z$ is {\bf surrounded} by $L$.

We define the distance between quasisimple loops $L_1$ and $L_1$
by $$d(L_1, L_2) = \inf ||\zeta_1 - \zeta_2||_\infty,$$ where the
infimum is taken over pairs $\zeta_1:S^1 \to \R^2$ and
$\zeta_2: S^1 \to \R^2$ of parameterizations of $L_1$ and $L_2$.

Let $D$ be a bounded domain.  Then the set $\Omega_D$ of all discrete
subsets of the set of quasisimple loops in $\overline D$ is a
metric space under the Hausdorff metric induced by the
metric $d(\cdot, \cdot)$ described above. Denote by $\mathcal F_D$
the Borel $\sigma$-algebra on $\Omega_D$.
%The laws of the \CLEk/ we
%construct below will be probability measures on $(\Omega_D, \mathcal F_D)$.
Most natural functions of $\Omega_D$ (such as the number of loops completely surrounding
a fixed disc, or the number of loops intersecting two disjoint open sets, or the event
that the outermost loop surrounding a fixed point is a simple loop) can be shown to be
$\mathcal F_D$ measurable.

The most natural definition of a random loop ensemble is a random variable whose
law is a probability measure on $(\Omega_D, \mathcal F_D)$.
Unfortunately, when defining \CLEk/ in arbitrary domains, it may not
be enough to consider quasisimple loops.  For example, consider the case that $D$
is a non-Jordan domain such as the square $(0,2) \times (0,2)$ minus the set $\{n^{-1} : n \in
\mathbb Z_+ \} \times (0,1)$.  When $\kappa = 8$, we expect \CLEk/ to be a
single space filling loop.  However, it is clear that any loop in the closure of $D$ which is
space filling in $D$ cannot be a continuous closed curve (since its $y$ coordinate
must oscillate between $0$ and $1$ infinitely many times).  Even when the original domain
is a Jordan domain and $\kappa < 8$, one may worry a priori that as we construct loops through an exploration
process we may create domains which are no longer Jordan domains.  It
will therefore be convenient to slightly expand our definition of loops.

A {\bf pinned loop} in $\H$ is a quasisimple loop $L$ in $\overline \H$ which intersects the origin at $0$
and has the property that if $\eta:S^1 \to \overline \H$ is a parameterization of the loop then
$\eta^{-1} \H$ is dense.  Roughly speaking, we now wish to define a ``conformal loop'' in $\C$
to be the image of a pinned loop under a conformal map $\phi$ from $\H$ to
another domain $D \subset \C$.  Note that the parameterization $\phi \circ \eta$ of this image
is well-defined on an open dense subset $U_0$ of $S^1$ but may not extend continuously to all of $S^1$.
If $U$ is the set of all points $s$ in $S^1$ for which $\phi \circ \eta$ can be extended to a
continuous function in a neighborhood of $s$, then it is clear that $U$ is the largest open set
to which $\phi \circ \eta$ can be continuously extended.  Formally, a
{\bf conformal loop} is a map $\zeta$ from a dense open subset $U$ of
$S^1$ into $\C$ (modulo monotone reparameterization) such that
\begin{enumerate}
\item There exists a conformal map $\psi$ from $\H$ to a superset of $\zeta(U)$
with the property that $\psi^{-1} \circ \zeta$ is the restriction
to $U$ of a function $\eta \in \overline{\mathcal C}$ with
$\eta^{-1} \H \subset U$.
\item The set $U$ is maximal, i.e., that there exists
no open proper superset $U'$ of $U$ such that $\zeta$ can be extended to a
continuous function of $U'$.
\end{enumerate}
%A {\bf conformal loop} is an equivalence class of such
%parameterizations, where $(\eta_1, U_1) \equiv (\eta_2, U_2)$ if
%$\eta_1 = \eta_2 \circ \psi$ for some topological automorphism
%$\psi$ of $S^1$ and $\psi(U_1) = U_2$.

If $L=(\zeta, U)$ is a conformal loop, we will sometimes abuse
notation slightly and use $L$ to denote the corresponding set $\overline{
\zeta(U)} \subset \C$ of points on the loop.  A point in $z \in \C
\setminus L$ is {\bf surrounded} by $L$ if $g(z)$ is surrounded by
$\eta$, where $\eta$ is as described above.  The reader may check that
this definition is the same for every $g$, and that the set of $z$
surrounded by $L$ is a union of bounded connected components of $\C
\setminus \L$.

Our initial approach
to defining \CLEk/ will be to define a coupling of \SLEbr \kappa-6/ processes called
a continuum exploration tree and to show that this process almost surely determines
a countable collection of {\it conformal loops}.  We will then show in Section \ref{s.confsymmetry}
that---if Conjecture \ref{c.reversible} holds---these loops are almost surely quasisimple and
the laws of the \CLEk/ may be equivalently described as measures on $(\Omega_D, \mathcal F_D)$.

%We define the distance between loops by $$d(L_1, L_2) = \inf
%||1_{U_1 \cap U_2}(\zeta_1 - \zeta_2)||_\infty,$$ where the $\inf$ is
%taken over all pairs of parameterizations $(U_1, \zeta_1)$ and $(U_2, \zeta_2)$.

\subsection{Continuum exploration trees} \label{s.contextree}

Fix initial values $W_0 = O_0 = 1 \in \partial \D$. Then Proposition
\ref{p.invariance} implies that an \SLEbr \kappa-6/ targeted at $a_1
\in \D$ and an \SLEbr \kappa-6/ targeted at another point $a_2 \in \D$
can be coupled in such a way that the corresponding growth processes
$K_t^{a_1}$ and $K_t^{a_2}$ agree (after a time change) up to the first time $t$
at which $a_1$ and $a_2$ are {\bf separated} (i.e., $a_1 \in K_t^{a_2}$ and $a_2 \in K_t^{a_1}$)
and evolve independently of one another after that time.
In other words, we may construct this coupling by first sampling an
\SLEbr \kappa-6/ process $K_t^{a_1}$ targeted at $a_1$, and then sampling an
\SLEbr \kappa-6/ process $K_t^{a_2}$ targeted at $a_2$ conditioned on it agreeing
with the first sample (up to time change) until the first time that $a_1$ and
$a_2$ are separated.

(We remark that if it were known that radial \SLEbr \kappa-6/
were a continuous path $\gamma$, then the separation time would
be the smallest $t$ for which $a_1$ and $a_2$ lie in distinct components of $\D \setminus \gamma([0,t])$.
We use instead the language of growth processes because
we have not proved that radial \SLEbr \kappa-6/ is a continuous path;
however, if we had a proof of Conjecture \ref{c.reversible},
then the continuity of radial \SLEbr \kappa-6/ would be an easy consequence,
and the above construction would describe a random path that ``branches'' at the point $\gamma(t)$.)

We can define a similar coupling inductively for sets $a_1, \ldots, a_k$ with $k > 2$ as follows.
To sample from such a coupling, first we choose the growth processes $K_t^{a_j}$ for all $j < k$
and then conditioned on these processes, we choose $K_t^{a_k}$ conditioned on having it agree with
each $K_t^{a_j}$ (after a time change) up until the first time $a_k$ is separated from $a_j$.
We can interpret this as an \SLEbkr/ process that ``branches'' at each
of the finitely many times that a pair $a_i$ and $a_j$ becomes separated for the first time.

In fact, we may repeat this procedure infinitely many times --- for some
countable dense sequence $a_1, a_2, a_3, \ldots$ of points in $\D$ ---
to obtain a coupling of \SLEbr \kappa - 6/ growth processes $K^{a_i}$ from
$1$ to $a_i$.  (Formally, if $\Omega^{a_i}$ is the space of continuous
driving parameters $W$ for growth processes targeted at $a_i$ and $\mathcal F^{a_i}$ is the
smallest $\sigma$-algebra which makes $W_t$ measurable for each fixed $t$, then this
coupling is a random variable in $\prod \Omega^{a_i}$ whose law is measurable with respect
to the product $\sigma$-algebra generated by the $\mathcal F^{a_i}$.)  Note that for every $i \not = j$, the
$K^{a_i}_t$ and $K^{a_j}_t$ agree almost surely (after a time change) until
the first time $t$ that $a_i$ and $a_j$ are {\bf separated} (i.e., $a_i \in K^{a_j}_t$
and $a_j \in K^{a_i}_t$), after which they
evolve independently.  We may view this collection of $K^{a_j}_t$ as a single random growth process
which branches at each of the countably many times that some $a_i$ becomes separated from
some $a_j$ for the first time.

Given a family of growth process $K^{a_i}_t$ (one for each $i \geq 1$) chosen from such a coupling,
we may almost surely uniquely (up to time change) define, for each point $z \in \overline \D$,
a growth process $K^z_t$ such that for each $a_i$,
the processes $K^{a_i}_t$ and $K^z_t$ agree (after a time change) until the first time
$a_i$ and $z$ are separated.  It is not hard to see that the joint law of the
processes $K^z_t$ (now defined for all $z \in \overline \D$) is independent of our
choice of countable dense set $\{a_i\}$.

The complete collection of processes $K^z_t$, for $z \in \overline \D$, is
called {\bf branching \SLEr \kappa - 6/}, or the {\bf continuum exploration tree}.

\subsection{Constructing loops from exploration trees} \label{CLEconstruction}

The $K^z_t$ are analogous to the exploration paths $T_v(A)$
that we defined in Section \ref{explorationtreesection} in the discrete setting.
(Recall Figure \ref{explorationtreefigure}.)  This section will describe
an algorithm for constructing a family of conformal loops from a continuum
exploration tree.  The algorithm is motivated by the discrete picture, in particular
Section \ref{s.excursion}.

Recall that in the discrete setting, the paths $T_v(A)$ trace part
(but not necessarily all) of each loop that surrounds $v$.  The reader
may observe (recalling Section \ref{s.excursion}) that if $T_v(A)$ begins to
trace such a loop at a vertex $w$ at step $m$,
then it will continue to trace that loop until the first time $n$ that\
$v$ and $w$ are separated by the path drawn thus far (i.e., $w$ fails to lie on the
boundary of $G_n$).  The portion of $T_v(A)$
between $m$ and $n$ is an arc of a loop.  The $m$ and $n$ are successive renewal times (as defined
in Section \ref{s.excursion}), and the corresponding graphs $G_m$ and $G_n$ have
monochromatic boundaries of opposite colors (as defined in Section \ref{s.excursion}).  We will now
make the analogous construction in the continuum setting.

Fix $z \in \D$ (we may assume $z = 0$, applying a
conformal map otherwise).  Write $\theta_t = \arg W_t - \arg
O_t$, where the branch of $\arg$ is chosen so that $\theta_t$ is continuous in $t$ and
$\theta_0 = 0$.  Since $\rho = \kappa-6$, using the definition for \SLEkr/ given in
\ref{s.radialSLEkr}, we see that the difference $\theta_t = \arg W_t - \arg
O_t$ is a real-valued process that evolves according to the diffusion
\begin{equation} \label{arclengthSDE} d\theta_t = \frac{\kappa-4}{2} \cot(\theta_t/2)\,dt + \sqrt{\kappa}\, dB_t
\end{equation}
in between those times $t$ for which $\theta_t$ is an integer
multiple of $2 \pi$.

We will now construct a sequence $L^z_j$ of nested conformal loops
surrounding $z$, such that each $L^z_{j+1}$ is surrounded by
$L^z_j$.  First, we call a time $t$ a {\bf loop closure time} if it
is the first time that $\theta_t$ hits a particular integer multiple
of $2\pi$ after the last previous time $s$ that it hit a different
multiple.  (Having $\theta_t$ change from an odd to an even multiple of
$2\pi$ corresponds to having the monochromatic boundary of $G_m$
and $G_n$ be of opposite colors in the discrete setting, as discussed above.)
Denote by $t^z_j$ the $j$th such time and by $s^z_j$ the
corresponding value of $s$ (which correspond to $m$ and $n$ in the discrete setting).
Denote by $A_j^z$ the component of $\D \backslash K_t$ that contains $z$.

We will see that the process $K_t$ between times $s^z_j$ and $t^z_j$ traces part of a conformal loop.
To begin to describe that loop, let $\phi: \D \setminus K_{s^z_j} \to \H$
be the composition of $g_{s^z_j}$ with a conformal map from $\D$ to $\H$ that sends $W_{s^z_j}$ to zero.

We now claim that for all
$t \in [s^z_j,t^z_j]$, the set $\H \setminus \phi K_t$ can be written, almost
surely, as the unbounded component of $\zeta([s^z_j,t])$
where $\zeta:[s^z_j, t^z_j] \to \overline \H$
is a continuous path which extends continuously to its endpoints. To see this, note
that the corresponding \SLEk/ in a Jordan domain is almost surely continuous by
\cite{\RohdeSchramm}.  For each fixed $s \in (s_j, t_j)$, the law
of the evolution of $K_t$ after time $s$ is that of an \SLEk/.  If $h_s$ is
a conformal map from $\H \setminus \phi K_s$ to $\H$ (with the hydrodynamic
normalization at infinity) then this implies that the growth of $h_s \phi K_t$
is indeed given by a continuous function $\zeta_s$ on $[s,t_j]$.  We can then set
$\zeta(t) = h_s^{-1} (\zeta_s(t))$ whenever $h_s^{-1}(\zeta_s(t))$ is well defined.
Since the $h_s$ converge uniformly to the identity, the claim follows from the fact
that the uniform limit of continuous functions is continuous.  We denote by $\gamma^z_j$
the path $\zeta$ described above for a fixed choice of $z$ and $j$.

\begin{figure}[t]\label{cleloopfig}
%\epsfbox[-45 43 100 240]
\begin{center}
\scalebox{.5}{\includegraphics{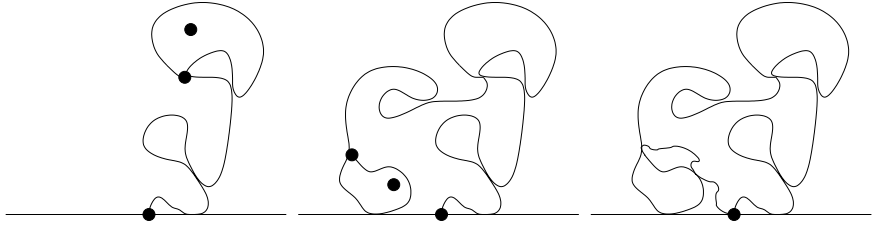}}
\caption {The left figure is a schematic drawing of $\gamma^z_j$ when $4 < \kappa < 8$.  The black dot
which does not lie on the path is $\phi(z)$.  The continuous path $\gamma^z_j$ begins
at $0$ (at time $s^z_j$) and ends when it reaches the first point $\gamma^z_j(t^z_j)$ (also shown as a black dot),
where $t^z_j$ is the first time that $\phi(z)$ is separated from the origin.  (For discrete
intuition, imagine that $\R$ is made of white hexagons, and the
path $\gamma^z_j$ has black hexagons on its left side and white hexagons on its right side.)
Conditioned on
this $\gamma^z_j$, the law of the remainder of the loop $\phi(L^z_j)$ is that of an \SLEk/ in $\H \setminus
\gamma^z_j([s^z_j,t^z_j])$ started at $\gamma^z_j(t^z_j)$ and ending at $0$.  The
middle figure is a schematic drawing of $\gamma^{z'}_j$ for a different value $z'$ (the black dot
which does not lie on the path); in this case, $\gamma^{z'}_j$ is an extension of the path $\gamma^z_j$.
The figure on the right indicates the union of all such extensions.  This is the complete loop
$\phi(L^z_j)$, a quasisimple closed loop in $\overline \H$ which begins and ends at $0$.}
\end{center}
\end{figure}

When $\kappa \leq 4$, we know that chordal \SLEk/ is a simple path which extends
continuously to its endpoints and does intersect the boundary of the domain
it is defined on except at these endpoints; from this we may conclude that the
second endpoint of $\gamma^z_j$ the same as the first endpoint, so that $\gamma^z_j$
is in fact a simple closed curve.  We then define $L^z_j$ to be the image
under $\phi^{-1}$ of $\gamma^z_j$.  Then $L^z_j$ is a conformal loop by definition.

When $\kappa > 4$, the path $\gamma^z_j$ stops the first time that
$z$ fails to lie on the component of $\H \setminus \gamma^z_j$ that
contains the first endpoint of $\gamma^z_j$ on its boundary.  We do
not expect $\gamma^z_j$ to be a closed loop.  However, we may
consider some $z' \in \D \setminus \cup_{s' < s} K_{s'}^z$.  Then it
is not hard to see that either $\gamma^z_j$ starts at a different
time from $\gamma^{z'}_j$ or $\gamma^z_j$ is (up to time change) a
proper sub-arc of $\gamma^{z'}_j$; moreover, the union of the
$\gamma^{z'}_j$ for which the latter holds is a quasisimple closed
loop. Given $\gamma^z_j$, its law is given by a chordal \SLEk/ (in
the appropriate component of $\H \setminus \gamma^z_j$) from the
last endpoint of $\gamma^z_j$ to the first endpoint of $\gamma^z_j$.
Now we define $L^z_j$ to be the image under $\phi^{-1}$ of this
loop. Again, $L^z_j$ is a conformal loop by definition.

To get a more intuitive understanding of the relationship between loops and trees,
recall that in the discrete setting, at each branch point that lies on a loop we have a notion
of a ``proper'' branch (which continues to follow the loop) and an ``improper'' branch (which
does not continue to follow the loop).  In the continuum setting, we can define branch points
analogously; given a pair $a_i \not = a_j$, and given $K_t^{a_i}$ and $K_t^{a_j}$ time changed so that they
agree up until the first time (call it $T$) for which $a_i$ and $a_j$ are separated, we say
that $K_T$ is the {\bf branch set} of $a_i$ and $a_j$, and $g_T^{-1} W_t$ (when it exists---i.e.,
when $g_T^{-1}$ extends continuously $W_t$)
is their {\bf branch point}.  When $\kappa \leq 4$, it is clear that the exploration tree cannot
branch in the middle of tracing a loop (since \SLEk/ is simple in this case).  However, when
$\kappa > 4$, the tree can branch while tracing the boundary of a loop; for each of the countably
many branch points of this form, there is a proper branch (which continues to trace the boundary of the loop)
and an improper branch (which turns into a region whose boundary is part of the loop).

When $\kappa \in (8/3, 8)$, $\beta \in [-1,1]$, and $\mu \in \R$ (with $\beta = 0$ if $\kappa=4$
and $\mu = 0$ otherwise), we define \CLEbk/ to be the collection of loops
of the form $L^z_j$ described above.  We will show in Section \ref{s.initialpointinvariance} that if
Conjecture \ref{c.reversible}
holds and $\kappa \in (4,8)$, then the law of the collection
of loops \CLEbk/ is independent of $\beta$ and $\mu$ --- and can thus be denoted \CLEk/
(see Proposition \ref{s.initialpointinvariance}).
(A similar fact for $\kappa \leq 4$ will appear in \cite{\SheffieldWerner}.)
In the absence of a complete proof of this fact, we will (somewhat arbitrarily)
declare \CLEk/ to be the ensemble of loops corresponding to $\mu = 0$ and either
$\beta = 1$ if $\kappa \in (4,8)$, or $\beta = 0$ if $\kappa \in (8/3, 4]$.
The following is immediate from Proposition \ref{p.invariance}
(and does not rely on Conjecture \ref{c.reversible}):

\begin{proposition} \label{p.clebkpartiallyconfinv}
Fix $\kappa \in (8/3, 8)$, $\beta \in [-1,1]$, and $\mu \in \R$ (with $\beta = 0$ if $\kappa=4$
and $\mu = 0$ otherwise).  Then the law of \CLEbk/ in a simply connected domain $D$ is
invariant under conformal automorphisms
of $D$ that fix the starting point of the exploration tree.
\end{proposition}

% The random collection
%$L^j_z$, defined for all $j \geq 1$ and all $z$ in some countable
%dense set, is called \CLEk/. When $\kappa \leq 4$, we may
%equivalently define the $L^j_z$ to be the boundaries of the open
%sets $A^j_z$.

\subsection{\CLE/ gasket and conformal radius distribution} \label{s.gasket}

The closure $\Gamma$ of the union $\cup L^z_1$ of all
``outermost'' loops is called the {\bf \CLE/ gasket}.  It is a
random closed set. It is not hard to see that conditioned on
$\Gamma$, the law of the non-outermost loops is given by an
independent \CLE/ in each component of the complement of $\Gamma$.
When $\kappa \leq 4$, the loops are the boundaries of the gasket.
Various properties of $\Gamma$ will be investigated in
\cite{\SchrammSheffieldWilson}.  As a prelude to
\cite{\SchrammSheffieldWilson}, we describe one such
property here.

Recall that when $\rho = \kappa-6$, the difference $\theta_t = \log
W_t - \log O_t$, in radial coordinates, evolves according to the
diffusion (\ref{arclengthSDE}) in between those times $t$ for which
$\theta_t= 0$ or $\theta_t = 2 \pi$. Let $\D$ be the unit disc and
let $A$ be the component of the complement of the gasket that
contains the origin---i.e., $A = \D \setminus K^z_t$, where $t =
t^z_1$. The {\bf conformal radius of $A$} (viewed from the origin)
is defined to be $|f'(0)|^{-1}$ where $f$ is any conformal map from
$A$ to $\D$ that fixes the origin.  Let $T$ be $-1$ times the log of
the conformal radius of $A$.  By the construction using branching
\SLEr \kappa-6/, it is clear that $T$ is equal to the time in a
\SLEr \kappa-6/ evolution that $\theta_t$ first hits $\pm 2 \pi$
when $\theta_0= 0$.

The probability density for the law of $T$ will be explicitly
computed in \cite{\SchrammSheffieldWilson}.  For now, we offer the
following:

\begin{proposition} \label{gasketconfrad}
There is a unique adapted (to filtration generated by the Brownian motion $B_t$)
and almost surely continuous random process $R_t$
on the interval $[0,2\pi]$ that is instantaneously reflecting at its endpoints (i.e., the total
amount of time spent at $0$ or $2\pi$ almost surely has Lebesgue
measure zero), satisfies $R_0 = 0$, and---in between times that
$R_t$ hits the boundary---evolves according to the SDE
(\ref{arclengthSDE}). The law of $T$ is equivalent to the law of
$\inf \{t : R_t = 2\pi \}$.
\end{proposition}

\begin{proof}
To prove the existence of a process $R_t$ with the properties
described above, we will show that $|\theta_t|$ (where $\theta_t$ is as defined
from the appropriate radial \SLEbkr/ as in Section \ref{CLEconstruction}),
has these properties.  We already observed in Section \ref{CLEconstruction} that the evolution of $|\theta_t|$ in
between times $[0,2\pi]$ is given by (\ref{arclengthSDE})); the fact that $|\theta_t|$ is
adapted and instantaneously reflecting can be seen by changing
coordinates to the setting of chordal \SLEkr/ in the half plane (Proposition \ref{p.radialchordalcoordchange})
and then recalling that Bessel processes of dimension $\delta > 1$ are instantaneously reflecting
(Proposition \ref{p.zerolebesgue}).

Next, if $R_t$ is any process with these properties, we may approximate it by
the process $R^\epsilon_t$ that evolves according to the SDE (\ref{arclengthSDE}) except that it
immediately jumps to $\epsilon$ each time the process hits zero.
Clearly, for each fixed $a \in (0, 2\pi]$, the law of $\inf \{t : R_t^\epsilon = a \}$ is
stochastically decreasing in $\epsilon$ (for $\epsilon < a$) and converges to the law of
$\inf \{t : R_t = a \}$ as $\epsilon \to 0$ (since $R_t$ is
instantaneously reflecting).  From this it is not hard to show that the $R_t^\epsilon$
converge in law (with respect to uniform topology on compact intervals) to $R_t$, and hence
the law of $R_t$ is uniquely determined.  The fact that the law of $T$ is equivalent to the law of
$\inf \{t : R_t = 2\pi \}$ is now immediate from the way that the loops were defined
(in terms of $\theta_t$) in Section \ref{CLEconstruction}.
\end{proof}

%The following somewhat technical fact will be useful in determining
%the law of $T$ explicitly:
%\begin{proposition}
%\end{proposition}

\subsection{Limiting cases: $\kappa = 8/3$ and $\kappa=8$}

As $\kappa \to 8/3$ from above, the dimension $\delta$ of
the corresponding Bessel $X_t$ tends to zero; the process
$X_t$ itself converges weakly to zero (with respect to the uniform topology
on compact intervals of time) as $\delta \to 0$.

From this and results of Section \ref{s.gasket}, it is not hard to see
that as $\kappa \to 8/3$ from above and $z \in D$ is fixed,
the conformal radius of the interior of each $L^z_j$, viewed from $z$, tends to zero in probability.  We thus
define  \CLEkk {8/3}/ to be the loop ensemble which
almost surely contains no loops.  (Alternatively, we could define every
point in $D$ to be its own loop.)

We remark that in the chordal \SLE/ setting (recall the definitions
in Section \ref{s.chordalsle}), as $\kappa$ tends to $8/3$ from
above, the $X_t$ becomes more tightly concentrated around zero, so
that in the limit we have $$W_t = \frac{4}{(2\rho+2)/\kappa)}B_t =
2\sqrt{\kappa}/(\kappa-4)B_t = \sqrt 6 B_t,$$ and hence, the driving
parameter of the \SLEkk \kappa, \kappa -6/ process converges weakly
to that of \SLEkk 6/ as $\kappa$ tends to $8/3$.  It is somewhat
natural that (as the loops become very small), the law of the
exploration path should converge to a process with the locality
property that \SLEkk 6/ has (given that renewal times --- i.e.,
times when $X_t = 0$ --- come with increasing frequency as $\kappa
\to 8/3$).

Next, we observe that when $\kappa \geq 8$ and $D$ is a Jordan domain, the
process \SLEr \kappa-6/ may be viewed as a space-filling
path $\gamma$ that starts and ends at the origin.  To see this, note that in the half
plane formulation we may first condition on
$W_t$ and $O_t$ up until a small stopping time $s$ for which $W_s \not = O_s$ almost surely.
Conditioned on $W_t$ and $O_t$ up until this time, the law of the remainder of the process
is given by an \SLEk/ (by Proposition \ref{p.chordalinvariance}), which is almost surely
space filling when $\kappa > 8$ \cite{\RohdeSchramm}.  Since we may take $s$ as small as we like,
it follows from \cite{\RohdeSchramm} that the process is a continuous path (when $D$ is a
Jordan domain) starting and ending at the same point.

We can therefore define \CLEk/ to be this single space-filling loop when $\kappa = 8$.
This loop can be approximated by choosing
two points $a$ and $b$ very close together on the boundary of $\H$ and drawing an \SLEk/ from $a$ to $b$.
It follows from \cite{\LSWUST} that \CLEk/ is the scaling limit of the path that
traces the boundary of the free uniform spanning tree (see \cite{\LSWUST} for a precise definition
of this paths).

In fact, we could define \CLEk/ analogously as a single space-filling loop for any $\kappa > 8$.
When $\kappa > 8$, however, we do not expect the law of this loop to be invariant under all
conformal and anticonformal automorphisms of $D$ (even though it is invariant under the
conformal automorphisms of $D$ that fix the starting point).  It is clear that
the law of \CLEk/ is invariant under all conformal and anticonformal automorphisms of $D$ when $\kappa = 8$ (because of
the uniform spanning tree interpretation) and when $\kappa = 8/3$ (trivially).  We will
discuss analogous invariance questions for the case $8/3 < \kappa < 8$ in Section \ref{s.confsymmetry}.

\section{Symmetry and uniqueness when $4 < \kappa < 8$} \label{s.confsymmetry}

\subsection{Continuity and starting point invariance of \CLE/} \label{s.initialpointinvariance}

This section will derive some consequences of Conjecture \ref{c.reversible} that
apply when $4 < \kappa < 8$.  Analogous questions for $8/3 \leq \kappa \leq 4$ will be dealt
with in a subsequent paper (\cite{SheffieldWerner}).  First,
we define {\bf boundary branching \SLEr \kappa-6/} to be a coupling of chordal \SLEr
\kappa-6/ processes targeted at each point in a countable dense set
of boundary points of $\D$ (instead of radial \SLEr \kappa-6/ processes targeted
at interior points of $\D$).  We may view this as a subset of the full
branching \SLEr \kappa-6/ tree, in which branching is only allowed to occur at points
on the boundary of $\D$.  Recall that the discrete analog of this tree traced all
of the loops that hit the boundary of the hexagonal graph (Proposition \ref{p.discreteboundarytree}).

In the continuum, it is also not hard to see that each of the conformal loops
traced by boundary branching \SLEr \kappa-6/ intersects the boundary of $\D$.  Let $\chi$ be the closure of the set of
points on loops traced by this process.  This $\chi$ is a random closed subset of the
\CLEk/ gasket.  By construction, conditioned on $\chi$,
the law of the remaining loops is given by an independent \CLEk/ in each component of the complement of
$\chi$, with appropriate starting points.  We will use this fact to prove the following:

\begin{proposition} \label{p.reversibleimpliesinvariance}
Fix $4 < \kappa < 8$,  $\beta \in [-1,1]$, and $\mu = 0$.  Suppose that Conjecture \ref{c.reversible} holds for $\kappa$.
Then the loops in a \CLEbk/ are almost surely quasisimple loops whenever $D$ is
a Jordan domain.  Moreover, the law of \CLEbk/ is invariant under conformal and
anticonformal automorphisms of $D$ and this law is independent of $\beta$.
\end{proposition}

\begin{proof} It can
be seen from Proposition \ref{p.radialchordalcoordchange} that if
chordal \SLEr \kappa-6/ is a continuous path when it is defined in
a Jordan domain, then radial \SLEr \kappa-6/ is
also; it then follows that the loops in \CLEk/ are almost surely quasisimple, since in this case
$\partial (\D \setminus K_t^{a_j})$ can be traced by a continuous curve, for every $j$ and $t$,
almost surely, and hence
any conformal map from $\H$ to $\D \setminus K_t^{a_j}$ extends continuously to $\R$.  This
also means that the processes $K^z_t$ are the hulls corresponding to actual continuous paths $\gamma^z(t)$
almost surely (i.e., for each fixed $z$, $\D \setminus K^z_t$ is almost surely the component of $\D \setminus
\gamma^z([0,t])$ which contains $z$, for all $t \geq 0$).

Now, to prove starting point independence, fix distinct points $a$ and $b$ on the boundary of $\D$ and let
$\psi$ be an anti-conformal map (i.e., a composition of a
conformal map and a reflection) from $\D$ to $\D$ such that $\psi(a)
= b$ and $\psi(b) = a$.  It will clearly be enough to show that the
law of the \CLEk/ is invariant under such a map $\psi$, since the group
of conformal and anticonformal automorphisms of $D$ is generated by functions $\psi$
of this form.  In fact, it will suffice if we can show that the law of the set of loops
in $\chi$ (i.e., the loops traced by the boundary-branching subtree of the exploration tree)
have a law which is invariant under such a $\psi$, since the law of all of the loops may be generated
inductively (as discussed above) from the law of $\chi$.

Let $T$ be an exploration tree $T$ started at $a$.  We define a subtree $T^{a,b}$ of the
boundary-branching exploration tree which
includes only the exploration path from $a$ to $b$ (i.e., the process $K^b_t$, which
is a continuous path if Conjecture \ref{c.reversible} holds) together with all paths obtained
by proper branches off of $T^{a,b}$.  In other words, this is the smallest subtree of the
exploration tree which traces all of the loops which are partially traced by the exploration
path from $a$ to $b$. See Figure \ref{clesymmetryfig}. (In the discrete analog,
given a set $A$ and boundary vertices $a= v_0$ and $b=v$, this would be the smallest subtree
of $T(A)$ that contains $T_v(A)$ together with all but one edge of every loop that $T_v(A)$
intersects.)

First, we claim that we can couple two instances $T^{a,b}_1$ and $T^{a,b}_2$ of the
exploration-tree-valued random variable $T^{a,b}$ in such a way that the set $\A_1$ of
quasisimple loops traced by $T^{a,b}_1$ is the same as the set $\psi \A_2$, where $\A_2$ is
the set of loops traced by $T^{a,b}_2$, almost surely.  To construct this coupling, we first arrange so
that the exploration
path $\gamma$ from $a$ to $b$ in $T^{a,b}_2$ is the image under $\psi$ of the exploration path from $a$
to $b$ in $T^{a,b}_2$ (which we can do if \SLEr \kappa-6/ has time reversal symmetry, since
the law of $\gamma$ is that of \SLEr \kappa-6/).

\begin{figure}[t]\label{clesymmetryfig}
%\epsfbox[-45 43 100 240]
\begin{center}
\scalebox{1}{\includegraphics{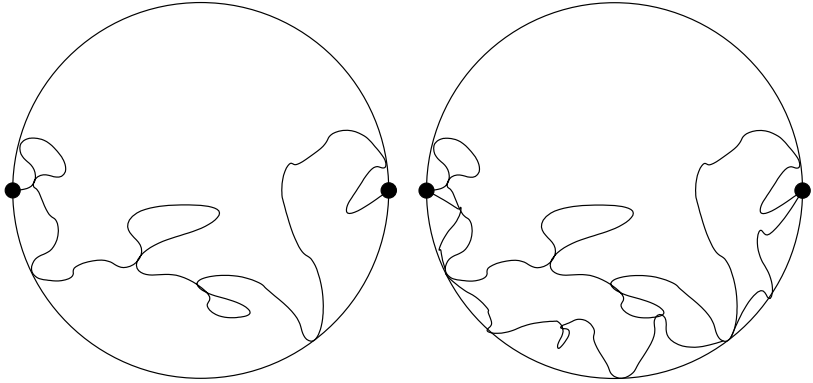}}
\caption {The left figure is a schematic diagram of an \SLEr \kappa-6/ $\gamma$ from the right boundary
point $a$ to the left boundary point $b$, where $4 < \kappa < 8$.  In this stylized drawing, the path hits the clockwise
arc from $a$ to $b$ at two points (instead of infinitely many).  These two points separate the
arc into three segments.  On the right, an \SLEk/ is drawn in each of the three regions of $\D \setminus \gamma$
that have these segments as boundaries.  The starting point for the \SLEk/ is the left endpoint of the segment
and the terminal point is the right endpoint.  This creates three quasisimple loops.  For discrete intuition, one
can imagine white hexagons along $\partial \D$ and the outside of each of the three loops and black hexagons
along the inside of each of three loops---so each component of the interior of the complement of the loops drawn
has a monochromatic boundary.  The exploration tree started at $a$ traces each of these loops in the
counterclockwise direction.}
\end{center}
\end{figure}

Next, let $D$ be a component of $\D \setminus
\gamma$ that includes a segment of $\partial D$ in its boundary;
let $y_1$ and $y_2$ be the first and last points on $\partial D$ hit
by the path $\gamma$; and let $\hat \gamma$ be the segment of $
\gamma$ which starts and ends at these points.  If $\kappa > 4$, then $\hat \gamma$ is
part of a loop of both $\phi \A_1$ and $\A_2$; in one of the trees $\psi T_1$ and $T_2$,
the remainder of the loop (obtained by following
the proper branches in the exploration process) is given by
an \SLEk/ from $y_1$ to $y_2$.  In the other tree, it is an \SLEk/ from $y_2$ to $y_1$ in
the same domain.  We couple the $T_i$ so that these two paths are equivalent up to time
parameterization (which we can do if \SLEk/ has time reversal symmetry).  Then
$\psi \A_1 = \A_2$ almost surely.

Now, let $D$ be a component of the interior of $\D$ minus the union of the loops in $\A_1$ such
that $D$ includes a segment of $\partial D$ in its boundary, and
let $y_1$ and $y_2$ be the first and last points of that segment $\partial D$.  Then the law of
the remainder of the boundary branching tree for $\L_2$ within the set $D$ is the image
under an anticonformal map $\psi_D$ from $D$ to itself (that maps $y_2$ to $y_1$ and $y_1$ to $y_2$)
of the law for $\L_1$.  We can then apply the same coupling we constructed above, where we replace
$\D$ with $D$ and the pair $(a,b)$ with $(y_1,y_2)$.  It is not hard to see that by repeating
this process for different segments of the boundary of $\D$, we obtain a coupling of the
boundary-branching exploration trees with the property that the loop ensembles $\B_1$ and $\phi \B_2$
generated by those trees agree almost surely.

%One may then construct a larger subset $\C_i$ of the loops in $\L_i$ by adding to
%$\B_i$ all of the loops generated by a new boundary branching exploration tree drawn in each
%of these components with an appropriate starting point.  By using a coupling of the type
%By repeating this process, one eventually obtains all of the loops in $\L_i$.
%By coupling the loop ensembles so that they agree at each stage of this process, we may arrange so that
%$\psi \L_1 = \L_2$ almost surely; hence, if $\L$ is a \CLEk/, then the law of $\psi \L$ is
%also that of a \CLEk/.
Finally, we need to show the independence of $\beta$.  Consider a
single exploration path $K^z_t$ of the exploration tree.  (Assume $z = 0$.)  Each
stopping time $T$ for which $O_T = W_T$ is a renewal time in the sense that conditioned on
$K^z_T$, the law of the loops in $\D \setminus K^z_t$ is given by a \CLEk/ in $\D \setminus K^z_t$.
Now, if we replace the exploration tree of this ensemble within $\D \setminus K^z_t$ by its
image under an anticonformal map, this will not change the law of the corresponding loop ensemble.

Let $\hat O_t$ denote the lifting of $\arg (W_t - O_t)$
to a continuous function on $\R$.  Recall (Proposition \ref{p.radialchordalcoordchange})
that radial \SLEr \kappa-6/ is defined so that if initial values $(W_0, \hat O_0)$ are
given where $\hat O_0 = 2 k \pi$ for some
integer $k$, and if $\psi$ is a conformal map (if $k$ is even) or an anti-conformal map
(if $k$ is odd) from $\D$ to $\H$ for which $\psi(W_0)=0$, then the image of
the corresponding radial Loewner evolution $K_t$ under $\psi$---up until
the first time $t$ that $\psi^{-1}(\infty)$ fails to lie on the boundary of
$\D \setminus K_t$---is given by
chordal \SLEbr \kappa-6/ with initial values $W_0=O_0=0$.

Fix $n > 0$.  One natural choice of stopping time $T = T_j$ (defined for each $j > 0$) is the first time
$t$ such the area of $D \setminus K^z_t$ is less than $j/n$ and $O_t = W_t$.  At each
such stopping time we may toss an independent coin with parameter $(\beta+1)/2$.  With
probability $(\beta+1)/2$, we continue the exploration tree with the original orientation (if $\hat O_0$
is an even multiple of $2\pi$) or the opposite orientation (if $\hat O_0$ is an odd multiple of $2\pi$).
An exploration path thus defined targeted at $z_1$ can be coupled with an exploration path thus
defined targeted at $z_2$ so that the two agree (after a time change) up until the first time that $z_1$ and $z_2$
are separated (since the definition of the stopping time makes no reference to choice of target point).
Coupling together paths of this form targeted at the dense set $\{a_j\}$, we obtain a variant
of the exploration tree which generates a set of loops with the same law as \CLEk/.
As $n \to \infty$, it is not hard to see that each exploration path in a tree thus defined converges
(in law with respect to uniform topology on compact time intervals) to \SLEbr \kappa-6/ (with $\mu = 0$).
Since for each $n$, the law of the loops is the same, the limit gives a coupling of the law of the loops
in a \CLEk/ with those of a branching \SLEbr \kappa-6/ process, which completes the proof.
\end{proof}

\subsection{Uniqueness result for boundary-intersecting \CLE/}

In this section we will focus on the case $4 < \kappa < 8$.
We will argue that if the scaling limit of an $O(n)$ models exists and
satisfies a few basic properties and conformal symmetries,
then it must be a \CLEk/ for some $\kappa > 4$, at least
when the model parameters are in the range for which the limiting loops hit the boundary.
In some sense we have already done this---namely, we proved Proposition \ref{p.Ktinvariance},
which suggests that, under reasonable conformal invariance hypotheses,
the scaling limit of the (non-oriented) exploration tree paths should
be \SLEr \kappa-6/---and hence the scaling limit of the exploration tree should be
branching \SLEr \kappa-6/.  Our definition of \CLEk/ was derived from that assumption.

Nonetheless, it is natural to wonder whether there are nice loop ensembles which do not naturally arise from
exploration trees in the same way.  Assuming Conjecture \ref{c.reversible}, we will give
a more complete axiomatic characterization of \CLEk/.  Recall that a loop is {\bf outermost}
if it is not surrounded by any other loop.

\begin{lemma} \label{givesallboundaryloops}
Fix $4 < \kappa < 8$ and suppose that Conjecture \ref{c.reversible} holds for that $\kappa$.  Let
$D$ be a Jordan domain.  Let the boundary branching exploration tree and \CLEk/ be coupled
as in the previous section.  Then the loops traced by the boundary branching exploration
tree are almost surely the only outermost loops in the \CLEk/ which intersect the boundary of $\D$.
\end{lemma}

%We also expect that there will almost surely be no non-outermost loops which intersect the boundary, but
%we will not need to prove this in order to prove Theorem \ref{t.cleuniqueness}.
The proof makes use of the following simple fact, which we will prove in Section \ref{s.approx}.

\begin{proposition} \label{p.epsilonBESconvergencelight}
Let $Z_T$ be number of times that a Bessel process $X_t$ of dimension $\delta > 1$
crosses the interval $[0,\epsilon]$ from bottom to top between time $0$ and time $T$.  Then
$\lim_{\epsilon \to 0} \epsilon \mathbb E Z_T = 0$.
\end{proposition}

\proofof{Lemma \ref{givesallboundaryloops}}
By Proposition \ref{p.reversibleimpliesinvariance}, we have (assuming that Conjecture \ref{c.reversible} holds)
that the law of \CLEk/ is invariant under the choice of starting point of the exploration tree.
At a renewal time of an exploration path (when $O_t = W_t$), the conditional law
of the loops in $\H \setminus K_t$ is that of a \CLEk/ in $\H \setminus K_t$.  We can therefore shift
the starting point of an exploration path at these renewal times without changing the
law of ensemble of loops generated at the end.  (We may think of this as exploring the same
loop ensemble beginning at a different point.)

We now define a process $(O_t', W_t')$ --- a invariant of our usual $(O_t, W_t)$ --- that
involves shifts of this form.  This process begins by evolving according to the diffusion
describing an \SLEkr/ process until time $t_1 = \inf \{t : -\epsilon \in K_t' \}$, where $K_t'$
is the Loewner hull generated by $(O'_t, W'_t)$. (Since we are assuming \SLEkr/ is continuous,
there is a continuous $\gamma$ corresponding to $K$ up to this point, and $t_1$ is the first time this path hits
$(-\infty, -\epsilon]$.)  At this point
the pair $(O'_t, W'_t)$ jumps $\epsilon$ units to the right (producing a discontinuity in $\gamma$).
Inductively, we let $t_{k+1}$ be the first time $t$ after time $t_k$
that $K_t$ hits $(-\infty, 0)$ (i.e., the corresponding continuous segment of $\gamma$
beginning at time $t$ hits a point in $(-\infty, 0)$). At each $t_k$, the process $(O'_t, W'_t)$ jumps $\epsilon$
units to the right and continues to evolve according to a the diffusion \SLEkr/ process until it jumps again
at $t_{k+1}$.  Thus, we can couple an $(O_t', W_t')$ with the driving parameters $(O_t, W_t)$ of
an ordinary \SLEkr/ process in such a way that $(O_t', W_t') = (O_t, W_t) + \epsilon(Y_t, Y_t)$ where $Y_t
= \inf \{k : t_k \geq t \}$.  Note also that $O_{t_k}' = W_{t_k}'$ for each $k \geq 1$ almost surely.

Let $Z_T$ be the number of times that $|W_t-O_t|$ hits $\epsilon$ for the first
time after the last time it hit $0$ before time $t = T$.  (In other words, $Z_T$ is the number
of upward crossings of the interval $[0, \epsilon]$ before time $T$.)  Now,
for each $t_k$, the probability $q$ that $|W'_t - O'_t|$ reaches the value $\epsilon$
before time $t_{k+1}$ is independent of $\epsilon$ and of the process $(W_t', O_t')$ up to
time $t_k$ (by scale invariance and the Markov property).  Thus, the expected number of $t_k$'s
occurring before time $T$ which are followed by such an upward crossing (at any point before $t_{k+1}$, which
may occur after time $T$) is exactly $q \mathbb E Y_T$.  We conclude that
$q \mathbb E Y_T \leq \mathbb E Z_T + 1$.  Proposition \ref{p.epsilonBESconvergencelight} implies
that $\lim_{\epsilon \to 0} \epsilon \mathbb E Z_T = 0$ and hence $\lim_{\epsilon \to 0} \epsilon \mathbb E Y_T = 0$.

Now, fix $a > 0$ and consider an \SLEkr/ stopped at $T = \inf \{t : -a \in K_t \}$.
%Clearly, $\inf \{t : -a' \in K_t \} \leq T$ whenever $a' < a + Y_T \epsilon$.
We may analogously define $T' = \int \{t : -a \in K'_t \}$.
Now consider a coupling of $(O_t, W_t)$ and the $(O_t', W_t')$ processes (one for each choice of
$\epsilon = 1/n$) in which the loops partially traced by
the corresponding Loewner evolutions belong to the same instance of \CLEk/.  Clearly, in such a coupling
the loops traced by $K_t$ up until time $T$ all touch the interval $(-a, 0)$.  However, the set $K'_{T'}$
contains {\em all} the loops that intersect $(-a, 0)$ (and possibly many more loops).
(To see this,
note that the corresponding $\gamma'$ hits the interval $(-a, 0)$ only finitely many times, and after
each such time there is a discontinuous jump to the right before the process starts again.)
The law of $T'$ therefore stochastically dominates
that of $T$ (since capacity is an increasing function of sets); since
$\epsilon Y_T \to 0$ in probability as $\epsilon \to 0$, it is not hard to see that $T' - T \to 0$ in probability
as $\epsilon \to 0$.

Now we claim that almost surely there is no loop in $\H \setminus K_T$ that intersects the interval
$[-a,0]$.  If there were such a loop with positive probability, then if we explored the same loop process
with the $(O_t', W_t')$ pair for some fixed $\epsilon$, as described above, it would have to
hit the loop before time $T'$ for all $\epsilon$, and thus $\lim_{\epsilon \to 0} T'$ would be
strictly larger than $T$.  Since $T' \geq T$ almost surely, and $T' \to T$ in law as $\epsilon \to 0$, we must
have $T' \to T$ almost surely as $\epsilon \to 0$.

Since each component of $\D \setminus \chi$ can be obtained as a union of $K_t$ processes of the form
described above, we see that there are almost surely no loops in any such component that intersect $\partial \D$.
\qed

For simplicity,
Theorem \ref{t.cleuniqueness}
will focus only on the law of the set $\mathcal L$ of outermost loops in a \CLEk/.
(Once this law is known, the law of the entire loop ensemble can be determined inductively.)

\begin{theorem} \label{t.cleuniqueness}
Suppose that $\L$ is a random countably infinite ensemble of quasisimple non-nested loops on $\H$
(formally, a probability measure on $(\Omega_D, \mathcal F_D)$, as defined in Section
\ref{s.initialpointinvariance}).
If Conjecture \ref{c.reversible} holds
for some $\kappa \in (4,8)$ and $\L$ is the set of outermost loops of a \CLEk/, then
$\L$ has the following properties:
\begin{enumerate}
\item \label{i.confinv} {\bf Conformal invariance:}
The law of $\L$ is invariant under conformal automorphisms of $\H$.
\item \label{i.boundint} {\bf Boundary intersection:}
The set $\L_\partial$ of loops of $\L$ that intersect $\R$ is almost surely non-empty,
and almost surely no loop in $\L_\partial$ hits any single point in $\R$ more than once.
%The loops of $\L$ are almost surely well separated.
%\item Given an arc $[a,b]$ of $\partial \D$, there is a almost surely a smallest closed interval
%$[c,d]$ of $\partial D$ for which there exists a loop $L$ that contains the points $c$
%and $d$ but does not intersect the set $[a,b]$.  Almost surely, there is only one
%such loop $L$ and each of $c$ and $d$ are hit only once by $L$.  Conditioned
%on the arc of $L$ from $c$ to $d$ that is closest to $[a,b]$, the law of the remain arc
%is \SLEk/ for some $\kappa > 0$.
\item \label{i.locfin} {\bf Local finiteness:} Given an interval $[a,b]$ and an open set $A \subset \H$ whose closure is
disjoint from $[a,b]$, there are almost surely only finitely many loops
which intersect both $A$ and the interval $[a,b]$.
\item \label{i.confmark} {\bf Conformal Markov property:}
Given $a$, $b$, and $A$ as in the previous item,
let $x$ be the right-most point in $[a,b]$ that lies in one of the (finitely many) loops $L$ that intersects both
$A$ and $[a,b]$.
Let $J$ be the counter-clockwise
arc of $L$ which begins at $x$ and ends at the first $y$ at which it hits $\partial A$.  Then $J$ almost surely does
not intersect $[a,x)$.  Given $J$ and the collection $\L_{(x,b]}$ of all loops that intersect $(x,b]$,
the conditional law of the counterclockwise arc of $L$ from $y$
to $x$ is given by an \SLEk/ from $y$ to $x$ in the component of
$$\D \setminus \overline{\cup \{L : L \in \L_{(x,b]}\} \cup J}$$
that has $[a,x]$ as part of its boundary.
\item \label{i.reprop} {\bf Renewal property:} Conditioned on the set $\L_{[a,b]}$ of loops of $\L$
intersecting an interval $[a,b]$ the law of the remaining loops in $\L$ is given by a product of independent random loop
ensembles in the (non-loop-surrounded) components of $\D \setminus \overline{\cup \L_{[a,b]}}$, each of which has the same
law as the original law of $\L$ conformally mapped to that component.
\end{enumerate}
Conversely (whether Conjecture \ref{c.reversible} holds or not), if $\L$
is any random countably infinite collection of quasisimple loops with the
properties listed above, then it must be a \CLEk/ for some $4 < \kappa < 8$.
\end{theorem}

Of course, if one wishes to avoid making explicit reference to \SLEk/ in
the conformal Markov property described above, one can replace the requirement that
the path from $y$ to $x$ is an \SLEk/ with the requirement that this random
path satisfies the hypotheses Theorem \ref{conformalmarkov}.

\begin{proof}

We begin by showing that if Conjecture \ref{c.reversible} holds for $\kappa$ with $4 < \kappa < 8$, then
\CLEk/ has the listed properties.  Proposition \ref{p.reversibleimpliesinvariance} implies \ref{i.confinv}, and
\ref{i.boundint} follows from easily the fact that \SLEk/ hits the boundary
and hits each point on the boundary at most once.  (The latter fact is an easy consequence of
continuity and time reversal symmetry of \SLEk/ and the fact that \SLEk/ hits each predetermined
boundary point with probability zero.  If $\gamma$ is a path from $0$ to $\infty$ in $\H$
that comes from a continuous Loewner evolution and $\gamma$ hits some $a \in \R \setminus \{0\}$
at two distinct times and hits another boundary point in $\R$ of the same sign in between these two times,
then the image of $\gamma$ under the inversion $z \to 1/\overline z$ cannot be a path
that comes from a continuous Loewner evolution.  Moreover, for each fixed rational $t$, the probability that
$\gamma$ hits the last boundary point that it hit before time $t$ for a second time after time $t$
is zero.  Hence, the probability that $\gamma$ hits a point in $\R$ at distinct times without
hitting another point in $\R$ in between those times is also zero.)

In the context of \ref{i.locfin},
the hypothesis that \SLEr \kappa-6/ is continuous implies that a chordal \SLEr \kappa-6/ from
$a$ to $b$ can have at most finitely many excursions away from the interval $[a,b]$ that intersect
$A$.  This together with Lemma \ref{givesallboundaryloops} implies \ref{i.locfin}.
Then \ref{i.confmark} and \ref{i.reprop} follow immediately from the conformal Markov
property and renewal properties of \SLEkr/ described in Proposition \ref {p.Ktinvariance}.

Now we proceed to the converse.  Suppose that $\L$ is a random ensemble of non-nested quasisimple
loops with all of the properties listed above.  Then by conformal invariance, the probability
that $\L$ contains a loop intersecting
an interval of $\partial \R$ is the same for all intervals.  Since this probability
must approach one as the length of the interval approaches $\infty$, we
conclude that $\L$ contains a loop intersecting each open interval of $\partial D$ with probability one.

Next, we will use a continuum analog of the construction of Propositions \ref{p.onechordalpathoftreebw} and
\ref{p.onechordalpathoftree} to construct a path from $0$ to
$\infty$ as follows.

Let $M_1, M_2, \ldots$ be an enumeration of the loops in $\L$ that intersect the
negative real axis.  For each $i\geq 1$, let $I_i$ be the interval
$$\left( \inf L_i \cap (-\infty,0), \sup L_i \cap (-\infty, 0) \right).$$
When the loop $L_i$ is counterclockwise oriented, let $A_i$ be
the portion of the loop $L_i$ that starts at $\sup
I_i$ and ends at $\inf I_i$. Suppose
that every $I_i$ is contained in some maximal $I_j$.  Then we may
consider the concatenation of the $A_i$, where $i$ ranges over those $i$ for which $I_i$ is maximal (i.e., for which
the interval $I_i$ is not contained in any distinct $I_j$).  More precisely, consider
any map from $(-\infty, 0]$ to $\H$ that maps each interval $I_i$ to the
corresponding loop segment $A_i$.  This is a map defined on an
open dense subset $\R$ to $\partial \H$, and local finiteness implies that it extends continuously
to a map from $[-\infty, 0]$ to $\partial \H$.  Let $Q_t$ be a parameterization of this path in
the opposite direction (so that $Q_0 = 0$).

We claim that the law of $Q_t$ must be that of \SLEr \kappa-6/ for
some $\kappa > 4$ when it is parameterized by capacity.
To prove this, it is enough to verify the hypotheses of
Proposition \ref{p.Ktinvariance}.

Let $t_a$ be the first time $t$ for which $Q_t \in (-\infty, a]$.  The renewal property implies that
conditioned on $Q_t$ up to time $t_a$, the law of the remainder of $Q_t$ is the same as the original
law of $Q_t$ after a conformal map $g_t:\H \setminus Q_t ([0,t_a]) \to \H$ that fixes $Q_{t_a}$ and
$\infty$.  In particular, the path $Q_t$, for $t \geq t_a$, remains in the closure of
the infinite component of $\H \setminus Q_t ([0,t_a]) \to \H$ almost surely.  Now, the conformal
Markov property of Proposition \ref{p.Ktinvariance} follows immediately from the conformal
Markov property cited here provided the stopping time is of the form $T = \inf \{t : Q_t \in A \}$,
where $A \subset \H$ is an open set whose closure does not contain $0$.  If $T'$ is any stopping
time such that almost surely $T' > T$ and $Q_{T'}$ lies on the same $A_i$ as $Q_T$ almost surely, then
it follows from Theorem \ref{conformalmarkov}.  The general result follows by noting that
for any stopping time $T''$ we can find a sequence of stopping times of the form $T'$ that converge to
$T''$ from below almost surely.

Finally, the conformal Markov properties implies that the Loewner evolution $W_t$
corresponding to $Q_t$ is continuous at all $t$ for which $Q_t \not \in \R$.  That this holds for
general $t$ follows easily from local finiteness.  We have now proved all the hypotheses of Proposition
\ref{p.Ktinvariance}.

By conformal invariance, the analogously defined $Q_t$---targeted at another point in $\R$
instead of $\infty$---will also have the law of \SLEr \kappa-6/ targeted at that point.  If
we consider a countable dense set of the boundary of $\overline \H$,
then we may define the union of the corresponding maps $Q$
targeted at these points to be a continuum exploration tree, which is a form
of branching \SLEr \kappa-6/; the fact that after the process branches the branches
evolve independently is immediate from the renewal property; it follows that
the law of the tree is the same
as the one given in Section \ref{CLEconstruction}.  The reader may now check
that the boundary-intersecting loops of $\L$ can be almost surely recovered from
this tree by applying the algorithm of Section \ref{CLEconstruction}.  By the renewal
property, the law of the boundary-intersecting loops determines the law of $\L$.

\end{proof}

\section{\SLEkr/ approximations and invariance}
\label{s.approx}

The authors in \cite{\SchrammWilson} proved a number of invariance
properties and coordinate changes for \SLEkr/ started with $O_0 \not
= W_0$ (so that $X_0 \not = 0$) and stopped at the first time $t$
such that $O_t = W_t$ (when $X_0 = 0$).  The main purpose of this section is to prove Proposition
\ref{p.radialchordalcoordchange}, which is essentially an extension of the analogous result in
\cite{\SchrammWilson} to times beyond the first time that $O_t = W_t$.

In order to prove this (without repeating all of the
calculations in \cite{\SchrammWilson}), we begin by
proving Proposition \ref{p.epsilonradialchordalcoordchange}, which shows
that \SLEbkr/ can be approximated by processes in which $O_t$
and $W_t$ are instantly pushed apart by a small fixed amount (which depends on $\mu$,
$\beta$, and $\kappa$) each time they collide. Proposition
\ref{p.epsilonradialchordalcoordchange} may also give the
reader more intuition about what the \SLEbkr/ processes are.  (Some of the
conjectures presented in Section \ref{s.openproblems} are based
on this intuition.)  Most of the
following exposition will focus on the case that $\mu = 0$, $\beta = 1$, and
$\kappa \not = 4$; the more general case will follow as a consequence of
this case.

\subsection{Approximate Bessel processes} \label{s.approxbessel}

Fix $\epsilon > 0 $.  Then we define an $\epsilon$-\BESdx/ process
$X^\epsilon_t$ to be a Markov process beginning at some initial
value $X^\epsilon_0 = x$ that evolves according to (\ref{BesselSDE})
except that each time it hits zero it immediately jumps to
$\epsilon$ and continues.  We may thus write
\begin{equation} \label{epsilonBesselSDE}
X_t^\epsilon = X^\epsilon_0 + \int_0^t
\frac{\delta-1}{2X^\epsilon_s}ds + B_t + J^\epsilon_t
\end{equation}
where $J^\epsilon_t$ is $\epsilon$ times the number of
$\epsilon$-jump discontinuities in $X_t^\epsilon$ up to and including time $t$.
(Note that if a jump occurs at $t$, then we write $X_t^\epsilon =
\epsilon$, so the process is upper semicontinuous.)  More generally,
a {\bf randomly jumping $\epsilon$-\BESdx/} process is a process in
which the jump sizes are random but the size of the jumps are almost
surely less than $\epsilon$, and the jump sizes are adapted to the
filtration generated by $B_t$. We denote the size of the $i$th jump
by $\epsilon_i$ and the time by $t_i$ and write
$$J^\epsilon_t = \sum_{t_i \leq t} \epsilon_i.$$
For simplicity, we also require that for each $t>0$ the set $\{i:
t_i \leq t \}$ is almost surely finite.  In particular, this
implies that the set of jump times $t_i$ is almost surely a discrete set.

\begin{proposition} \label{p.epsilonBESconvergence}
For each $\epsilon > 0$, let $X^\epsilon$ denote any randomly
jumping $\epsilon$-\BESdx/ process. As $\epsilon \to 0$, the
$X^\epsilon_t$ converge in law to a \BESdx/ with respect to the
$L^\infty$ metric on a fixed interval $[0,T]$, with $T>0$.
We also have almost surely, as $\epsilon \rightarrow 0$,
\begin{enumerate}
\item \label{geq1} $J_T^\epsilon \rightarrow 0$ if $\delta > 1$.
\item \label{eq1} $J_T^\epsilon \rightarrow l^0_t$ if $\delta = 1$.
\item \label{l1}$J_T^\epsilon \rightarrow \infty$ if $0 < \delta <
1$.
\item \label{l1squared} $J^{\epsilon^2}_T := \sum_{t_i \leq T} \epsilon_i^2 \to
0$ for all $\delta > 0$.
\end{enumerate}

\end{proposition}

\begin{proof}
We will deduce all of these results from the existence of a
continuous Bessel process $X_t$ of dimension $\delta$, which is a
strong solution to (\ref{BesselSDE}) away from zero, satisfies
Brownian scaling, and almost surely hits zero on a set of zero
Lebesgue measure (Propositions \ref{p.zerolebesgue} and
\ref{p.brownianscaling}).

First, we will construct a coupling of the process $X^\epsilon_t$
with a \BESdx/ process $X_t$ in such a way that that the two processes
agree when certain intervals of
time are excised from latter.  We will use $X$ use it to construct
$X^\epsilon$ as follows. First, set $X^\epsilon_t = X_t$ for $t \in
[0, t_1)$, where $t_1$ is the time at which $X_t$ first hits zero.
(Note that $t_1=0$ and the interval is empty if $x=0$.) At this
point we sample $\epsilon_1$ from the law of $\epsilon_1$ (in the
$X^\epsilon$ process) conditioned on our choice of $X^\epsilon_t$ up
to time $t_1$. Now we inductively define times $t_i$ and $s_i$ as
follows. Let $s_0 = 0$ and let $s_i$ be such that $t_i + s_i$ is
equal to the first time $t$ after $t_i + s_{i-1}$ for which $X_t =
\epsilon_i$. Then we define $X_{t}^\epsilon = X_{t+s_i}$ for $t \in
[t_i, t_{i+1})$, where $t_{i+1}$ is the first time $t>t_i$ for which
$X_{t+s_i}=0$. Then we choose $\epsilon_{i+1}$ from the law of
$\epsilon_{i+1}$ (in the $X^\epsilon$ process) conditioned on our
choice of $X^\epsilon_t$ up to time $t_{i+1}$ and continue.

We may think of $X^\epsilon_t$ as being obtained from $X_t$ by
``skipping'' the intervals of time $[t_i + s_{i-1}, t_i + s_i)$. On
each such interval, $X_t=0$ when $t$ is the left endpoint, $X_t =
\epsilon_i$ when $t$ is the right endpoint, and $0 \leq X_t <
\epsilon_i$ for other times $t$ in the interval.  Hence the skipped
time during finite interval $[0,T]$ is a subset of the set of times
$t$ for which $X_t \leq \epsilon$.  The total measure of the latter
set tends to zero as $\epsilon \to 0$, so it is clear that
$X^\epsilon$ and $X$ agree on $[0,T]$ up to translation of time by
an amount that tends to zero as $\epsilon \to 0$. Since $X$ is
almost surely continuous (and hence uniformly continuous on $[0,T']$
for any fixed $T'>T$), this implies that in the couplings above
$X^\epsilon \rightarrow X$ uniformly on $[0,T]$ almost surely as
$\epsilon \to 0$.  In particular, this implies that $X^\epsilon \to
X$ in law as claimed.

Now \ref{geq1} follows immediately from the uniform convergence of
$X^\epsilon$ to $X$ on $[0,T]$ together with convergence of the
other terms (besides $J^\epsilon_t$) on the right hand side of
(\ref{epsilonBesselSDE}) to the corresponding terms in
(\ref{BesselSDE}).  We get \ref{l1} by a similar argument and the
fact (from Proposition \ref{p.semimartingale}) that the integral on
the right hands side of (\ref{epsilonBesselSDE}) tends to $-\infty$
as $\epsilon \to 0$.  When $\delta = 1$, this integral is equal to
zero, so the \ref{eq1} follows from (\ref{epsilonBesselSDE}) and
(\ref{BesselSDE1}).

Next, let $t_a$ be the first time $t$ for which $J^{\epsilon^2}_t
\geq a$ for some fixed $a> 0$ and observe that \ref{l1squared} will
follow if we can show that $t_a \to \infty$ in probability, as
$\epsilon \to 0$, for each $a > 0$.  Let $A$ be a random variable
whose law is that of the first time that a \BESdx/ hits zero when
$x=1$. Let $m$ and $v$ be the mean and variance of $\max(A,1)$.
Write $b_i = \max(t_{i+1} - t_i, \epsilon_i^2) - m\epsilon_i^2$.
Given $\epsilon_i$, the law of $b_i$ is that of $\epsilon^2
(\max(A,1)-m)$ and has zero mean and variance $v \epsilon_i^4$. Thus
$$M_t = \sum_{i: t_i \leq t} b_i$$ is a martingale.

The variance of $M_{t_a}$ is $O(\epsilon^2)$.  Hence $\sum_{i: t_i
\leq t_a} \max(t_i - t_{i+1}, \epsilon_i^2) - \sum m\epsilon_i^2$
tends to zero in probability as $\epsilon \to 0$, which implies
$$\sum_{i: t_i \leq t_a} \max(t_i - t_{i+1}, \epsilon_i^2) \to ma$$
in probability.  The left hand side is bounded by $N_{t_a}$ where
$N_t$ is the Lebesgue measure of the set of times in $[0,t]$ that
are at most $\epsilon$ from a time $s$ for which $X_s = 0$.  Since
$N_t$ tends to zero in probability for each fixed $t$, as $\epsilon
\to 0$, we must have $t_a \to \infty$ in probability.
\end{proof}

\proofof{Proposition \ref{p.epsilonBESconvergencelight}}  This is immediate
from the fact that when $X_t$ and $X_t^\epsilon$ are coupled
as in the proof above, we have $J_t^\epsilon \geq \epsilon Z_t$ almost surely
for all $t$.
\qed

\subsection{Approximations to chordal \SLEkr/} \label{s.approxslekr}

\medskip

Let $O_t$ and $W_t$ be the parameters of an \SLEkr/ processes.
We will first consider the
case that $\rho < -2$, so that $\delta < 1$, and we assume also
that $\delta>0$. Take $X^\epsilon_t$
to be a randomly jumping $\epsilon$ \BESd/. Then we may write

\begin{eqnarray*}
X^\epsilon_t &=& \int_0^t \frac{\delta-1}{2X^\epsilon_s} ds + B_t +
J^\epsilon_t \\
Y^\epsilon_t &:=& \frac{2}{\delta-1}(X^\epsilon_t - B_t) = \int_0^t
(X^\epsilon_s)^{-1}ds + \frac{2}{\delta-1} J^\epsilon_t \\
\end{eqnarray*}

By Proposition \ref{p.epsilonBESconvergence} and (\ref{BesselSDEsmalldelta}),
the processes $X^\epsilon$ and $Y^\epsilon$ converge in law to $X$ and $Y$ as
$\epsilon \to 0$ (with respect to the $L^\infty$ metric on finite
intervals), where $X$ is a \BESd/ and $Y_t = \PV \int_0^t X^{-1}_s ds$.  This
implies that the following converge (in law with respect to $L^\infty$ metric
on compact intervals) to $O_t$ and $W_t$ as $\epsilon \to 0$:

\begin{eqnarray*}O^\epsilon_t & := & \frac{-2}{\sqrt \kappa}Y^\epsilon_t =
\frac{-2}{\sqrt \kappa} \int_0^t (X_s^\epsilon)^{-1} ds +
\frac{-4}{\sqrt \kappa (\delta - 1)} J^\epsilon_t \\
W^\epsilon_t &:=& O^\epsilon_t + \sqrt \kappa X^\epsilon_t =
\bigl(\frac{-2}{\sqrt \kappa} + \frac{\sqrt \kappa (\delta -
1)}{2}\bigr)\int_0^t (X^\epsilon_t)^{-1}ds + \\ & & \bigl(
\frac{-4}{\sqrt\kappa (\delta - 1)} + \sqrt \kappa \bigr)
J_t^\epsilon - \sqrt \kappa B_t \\
\end{eqnarray*}
Equivalently we may write
\begin{equation*}
\begin{pmatrix} O^\epsilon_t \\ W^\epsilon_t \\
\end{pmatrix} =
\begin{pmatrix} o_1 & o_2 & o_3 \\ w_1 & w_2 & w_3 \\ \end{pmatrix}
\begin{pmatrix} \int_0^t (X_s^\epsilon)^{-1} ds \\  J^\epsilon_t \\
B_t \\
\end{pmatrix},
\end{equation*}
where
\begin{equation*}
\begin{pmatrix} o_1 & o_2 & o_3 \\ w_1 & w_2 & w_3 \\ \end{pmatrix} =
\begin{pmatrix}
\frac{-2}{\sqrt \kappa} &
\frac{-4}{\sqrt \kappa (\delta - 1)} & 0 \\
\bigl(\frac{-2}{\sqrt \kappa} + \frac{\sqrt \kappa (\delta -
1)}{2}\bigr) &  \bigl( \frac{-4}{\sqrt\kappa (\delta - 1)} +
\sqrt \kappa \bigr) & - \sqrt \kappa \\
\end{pmatrix}.
\end{equation*}
When we substitute $\delta = 1 + \frac{2(\rho+2)}{\kappa}$ this
becomes
\begin{equation} \label{e.matrixforOandW}
\begin{pmatrix}
\frac{-2}{\sqrt \kappa} &
\frac{-2\sqrt \kappa}{ \rho+2}& 0 \\
\bigl(\frac{-2}{\sqrt \kappa} + \frac{\rho+2}{\sqrt \kappa}\bigr) &
\bigl( \frac{-2 \sqrt \kappa}{\rho+2} +
\sqrt \kappa \bigr) & - \sqrt \kappa \\
\end{pmatrix} =
\sqrt \kappa \begin{pmatrix} \frac{-2}{\kappa} &
\frac{-2}{\rho+2} & 0 \\
\frac{\rho}{\kappa} & \frac{\rho}{\rho+2} & - 1 \\
\end{pmatrix}.
\end{equation}

The values $O^\epsilon_s \leq W^\epsilon_s$ evolve as they would in
an ordinary \SLEkr/ process up until the first time $t$ for which
$W^\epsilon_t = O^\epsilon_t$.  At the $i$th time that this happens,
the value $W^\epsilon$ jumps by $\tilde \epsilon_i :=w_2 \epsilon_i
= \frac{\sqrt \kappa \rho}{\rho+2} \epsilon_i$ and $O^\epsilon$
jumps by the amount $\hat \epsilon_i : =o_2 \epsilon_i = -\frac{2}{\rho} \tilde
\epsilon_i$. The jump in $W^\epsilon_t$ corresponds to a
discontinuity in the path $\gamma$ that the Loewner evolution
describes.  Equivalently, a jump in $W^\epsilon_t$ (when time is
parameterized by capacity) corresponds to having $\gamma$ ``trace
the boundary of $K_t$'' until the image of its tip under $g_t$ has
moved $\tilde \epsilon_i$ units to the right; we also have $O^t$ move
along the boundary of $K_t$ until its image under $g_t$ has moved
$\hat \epsilon_i = \frac{-2}{\rho} \tilde \epsilon_i$ units to the
right.  (Since $\rho < -2$, we have $\hat\epsilon_i < \tilde \epsilon_i$.)

The above yields a construction of an $\epsilon$ approximation of an
\SLEkr/ process defined for all time in terms of \SLEkr/ processes
starting at $W_0 \not = O_0$ and stopped the first time $W_t = O_t$.
Observe that the ratio fixed $\hat \epsilon_i/\tilde \epsilon_i$ in the above
construction is canonical since $J^\epsilon_t \to \infty$ as
$\epsilon \to 0$ when $\rho < -2$, by Proposition
\ref{p.epsilonBESconvergence}.  If we fix any other ratio for
$\hat \epsilon_i/\tilde \epsilon_i$ and take $\tilde \epsilon_i \to 0$, then
$|W^\epsilon|$ will converge in law to $\infty$ instead of to a
continuous process $W$.

If $\rho \geq -2$ (so that $\delta \geq 1$), we define $O^t_\epsilon$ and $W^t_\epsilon$ as in
(\ref{e.matrixforOandW}) but with $w_2 = 0$ and $o_2 = -\sqrt
\kappa$.  Thus, as before, the values $O^\epsilon_s \leq
W^\epsilon_s$ evolve as they would in an ordinary \SLEkr/ process up
until the first time $t$ for which $W^\epsilon_t = O^\epsilon_t$.
However, at that point, the value $O^\epsilon$ jumps by $\hat
\epsilon := \sqrt \kappa \epsilon$ and $W^\epsilon$ does not jump.
Thus $W_t^\epsilon$ is almost surely continuous. Because, by
Proposition \ref{p.epsilonBESconvergence}, $J^\epsilon_t \to 0$ in
this setting, the process still converges to \SLEkr/ in the limit:

\begin{proposition}
In both the $0 < \delta < 1$ and the $\delta \geq 1$ settings
discussed above, the processes $W^\epsilon$ converge in law (with
respect to the $L^\infty$ metric on any compact interval $[0,T]$) to
the driving parameter of \SLEkr/ as $\epsilon \to 0$.
\end{proposition}

\subsection{Approximations to radial \SLEkk \kappa; \kappa - 6/}
We can define an approximation to radial \SLEr \kappa-6/ the same
way we did before in the chordal case in Section \ref{s.approxslekr}
except that the jumps are along the unit circle instead of the real
line. That is, if $\kappa < 4$ (so that $\rho < -2$), then at the
$i$th time $W^\epsilon_t = O^\epsilon_t$ the value $\arg W^\epsilon$
(instead of the value $W^\epsilon$) jumps by $ \tilde \epsilon_i
:=\mathcal E w_2 \epsilon_i$ and $\arg O^\epsilon$ jumps by the
amount $\hat \epsilon_i := \mathcal E o_2 \epsilon_i$ where $\mathcal E \in \{-1,1\}$
is $1$ if $W_t$ collided with $O_t$ on the clockwise side of $O_t$ and
$-1$ if $W_t$ collided with $O_t$ on the counterclockwise side of $O_t$.
When $\kappa > 4$ (so that $\rho > -2$), the value
$O^\epsilon$ jumps by $\hat \epsilon := \mathcal E \sqrt \kappa
\epsilon$ units to the left and $W^\epsilon$ does not jump.

\begin{proposition} \label{p.epsilonradialchordalcoordchange}
Let $W^\epsilon$ and $O^\epsilon$ be the processes discussed above, where
initial values $a \in \partial \D$ for $W^\epsilon_0$ and $b \in
\partial \D$ for $O^\epsilon_0$ are given.  As $\epsilon \to 0$, the $W^\epsilon_t$
described above converges in law (with respect to the $L^\infty$ metric
on fixed compact intervals of time) to a random process $W_t$.  Moreover
if $\psi$ is a conformal map from $\D$ to $\H$, then the image $\tilde K_t$ of
the corresponding Loewner evolution $K_t$ under $\psi$ is (up to
a time change) a growing family of closed sets given by a Loewner
evolution whose driving parameter converges in law
(with respect to the $L^\infty$ norm on intervals of the form $[0,T]$,
where $T$ is any bounded stopping time that satisfies $T < \overline T
:= \inf \{t : \psi(0) \in \tilde K_t \}$ almost surely)
to the driving parameter of \SLEr \kappa-6/ with initial values $W_0 =
\psi(a)$ and $O_0 = \psi(b)$.
\end{proposition}

\begin{proof}
We will prove the latter statement first.  Let $\phi_0: \D \to \H$ be the conformal map given by
$\phi_0(z) = \frac{2i}{z+1} - i$.  This satisfies $\phi_0(1) = 0$.
Given $a,b \in \R$, define $\phi_{a,b}(z)$ by $\phi_{a,b}(z) = a\phi_0(e^{ib} z) - a \phi_0(e^{ib})$.
We also have $\phi_{a,b}(1) = 0$, and $a$ and $b$ parameterize the set of conformal maps from
$\D$ to $\H$ with this property.

Take $\tilde K_s = \phi K_s$.  In this proof we will use $s$ to denote the
time of the radial process and define $t =t(s)$ be the half-plane
capacity of $\tilde K_s \subset \H$.  (We also write $s(t)$ for the choice
of $s$ for which $t=t(s)$.) Denote by $g_s: \D \setminus
K_s \to \D$ the Loewner conformal map at time $s$ in the disc, and
by $\tilde g_t: \H \setminus \tilde K_{t(s)} \to \H$ the conformal Loewner map
in the half-plane with the hydrodynamic normalization (i.e,
$\lim_{|z|\to\infty} |g_t(z)-z| = 0$).

Write $\psi_t = \tilde g_t \phi g_{s(t)}^{-1}$.  Then $\psi_t$ is a
conformal map from $\D$ to $\H$. We can extend this map to the
boundaries and write $\tilde W^\epsilon_t = \psi_t W^\epsilon_t,$
and $\tilde O^\epsilon_t = \psi_t O^\epsilon_t.$
Write $f_t = g_{s(t)}/W_{s(t)}$ and $\tilde f_t = \tilde g_t - \tilde W_t$.
These maps are normalized to send the tip of $K_t$ to $1$ and $0$,
respectively.  Write $\overline \psi_t = \tilde f_t \phi f_t^{-1}$.

This $\overline \psi_t: \D \to \H$ is a map that evolves, as $t$ grows,
within the two parameter family of conformal maps from $\D$ to $\H$
that map $1$ to $0$.  It can be described by the pair $a=a_t$ and $b=b_t$
as described above.
%by the
%parameters $(\lambda_t, \eta_t)$ where $\lambda_t = \psi_t'(0)/i \in
%\R$ and $\eta_t = \overline \psi_t^{-1}(\infty) \in
%\partial \D$.
It is not hard to work out the SDE describing the time evolution of
$a_t$ and $b_t$ (which is similar to what is done in
\cite{\SchrammWilson}).  For our purposes, it is enough to observe
that they vary continuously in $t$ (which is immediate from the
fact that the $g_t$ vary continuously with $t$).  Fix positive real constants $c$, $C$, and $d < \pi$
and let $T_{c,d, C}$ be the first time that $\tilde W_t =
\tilde O_t$ after the time $\inf \{t: t=C \text{ or } a_t = c \text{ or } |b_t| = d \}$.

Let $F_t = F_{a_t,b_t} (\theta) = \phi_{a_t,b_t}(e^{i\theta})$ be the real-valued extension of $\phi_{a_t,b_t}$
to $\partial \D$ (parameterized by the real angle parameter $\theta$).
Clearly, each of the derivatives of $F_t$ at $\theta=0$ is uniformly bounded up to time $T_{c,d,C}$.

By Proposition \ref{p.slecoordinatechange1}, $W^\epsilon_t$ and
$O^\epsilon_t$ initially evolve (up to a time change) according to the rule of
an \SLEk/ from $W^\epsilon_0$ to $O^\epsilon_0$. However, whenever
$\tilde W_t$ and $\tilde O_t$ collide there are jumps in both
$\tilde W_t$ and $\tilde O_t$ of size $F_t(\tilde \epsilon_i)$
and $F_t(\hat \epsilon_i)$.

Because of the uniform bounds on the derivatives of $F_t$, the sizes
of the $i$th jumps of $W^\epsilon_t$ and $O^\epsilon_t$ are $\tilde
\epsilon'_i$ and $\hat \epsilon'_i$ (for $\epsilon'_i = F'_t(\epsilon_i)$) plus an error which is
$O(\epsilon^2)$. Letting $\epsilon$ tend to zero, the fact that convergence
holds on the interval $[0,T_{c,d,C}]$ follows immediately from Proposition \ref{p.epsilonBESconvergence}.
The latter proposition statement then follows from the fact that $a_t < \pi$ and $b_t < \infty$ up until time
$\overline T$, and hence given any bounded stopping time $T$ which is almost surely less than $\overline T$,
we can choose $c$, $C$, and $d$ large enough so that the
probability that $T>T_{c,D,C}$ is arbitrarily close to zero.

To prove the first statement in the proposition, we first note that the proof
above implies the convergence in law of $W^\epsilon_t$,
at least up to some positive stopping time $T$ for which $O^\epsilon_T = W^\epsilon_T$ almost surely, to some
limiting process (namely, the driving parameter of the $\phi$ pre-image
of chordal \SLEr \kappa-6/).  By the Markovian property of the pair $(W^\epsilon_t, O^\epsilon_t)$,
and the fact that these stopping times are renewal times, this convergence holds for a stopping
time whose law is an independent sum of $k$ stopping times of this form.  The first statement then follows by taking
$k \to \infty$ and noting that the probability that such a sum is less than any fixed constant
tends to zero in $k$.
\end{proof}

%When $\kappa \in (8/3,8) \setminus \{4\}$, we can formally define
%{\bf radial \SLEk/} to be the Loewner evolution driven by $W_t$
%where the law of $W_t$ is the limit of the of laws of the
%$W^\epsilon_t$ defined above as $\epsilon \to 0$, noting

We have now essentially proved Proposition \ref{p.radialchordalcoordchange}.

\proofof{Proposition \ref{p.radialchordalcoordchange}}
When $\kappa \not = 4$ and $\beta=1$ and $\mu=0$, this existence of the process is immediate from Proposition
\ref{p.epsilonradialchordalcoordchange}, and the uniqueness is trivial.
 When $\delta \not = 1$, and
$\beta$ is general, similar arguments to those in Section
\ref{s.approxslekr} can be used to approximate \SLEr \kappa-6/ with
randomly jumping processes, to define branching analogs of skew
\SLEr \kappa-6/, and to show that these processes are invariant
under Mobius transformations of the domain that fix the starting
point.  The only difference is that in the $\epsilon$
approximations, each time $W^\epsilon_t$ and $O^\epsilon_t$ collide, instead of always
adding the appropriate $\tilde \epsilon_i$ to $W^\epsilon_t$ and $\hat \epsilon_i$ to
$O^\epsilon_t$, we add these values with probability $\frac{1+\beta}{2}$ and
subtract them with probability $\frac{1-\beta}{2}$.

When $\kappa = 4$, $\beta = 0$, and $\mu \in \R$, we obtain laws of
the driving parameters for chordal \SLEbr \kappa-6/ as weak limits of the laws of the corresponding
processes for $\kappa \not= 4$ (note the convergence of the corresponding
L\'evy processes described in Section \ref{s.levy}).  We may thus obtain laws for
radial \SLEbr \kappa-6/ as the corresponding limits of
the laws for radial \SLEbr \kappa-6/ as $\kappa \to 4$.  In all cases, the uniqueness
is trivial once existence is shown.
\qed

We briefly remark that when $\kappa \not = 4$ and $\mu = \beta = 0$, it is natural to modify the
approximation of \SLEr \kappa-6/ so that instead of adding some
$w_2\epsilon_i$ to $W^\epsilon_t$ and $o_2\epsilon_i$ to $O^\epsilon_t$ at a jump time
we leave $W^\epsilon_t$ fixed and add $(o_2-w_2)\epsilon_i$ to $O^\epsilon_t$.  This amounts
to shifting the whole process by $\pm w_2 \epsilon_i$.  The variance
of the sum of these random shifts is the expected sum of the squares
of the $w_2 \epsilon_i$, and it tends to zero as $\epsilon \to 0$, by
Proposition \ref{p.epsilonBESconvergence}.  It is then natural to rescale
and replace $(o_2 - w_2)\epsilon_i$ with $\epsilon_i$ --- so that at each jump, $O^\epsilon_t$ jumps
by $\pm \epsilon_i$.  This definition makes sense when $\kappa = 4$ as well; in this case,
the interested reader may check that one can obtain \SLEbr \kappa-6/ for general
$\mu$ as a limit by replacing $(\pm \epsilon_i)$ with $(\mu \epsilon -\epsilon_i,
\mu \epsilon + \epsilon_i)$.

%\subsection{Proof of Lemma \ref{givesallboundaryloops}}

\section{Height functions and other lattices}
\label{s.otherdiscrete}

\subsection{Height functions: continuity and monotonicity}

We now return to the discrete setting of Section \ref{s.discretesection}.
That is, we let $G$ be a hexagon graph with a set $F$ of hexagonal faces and a fixed vertex $v_0$ on its boundary
and a directed edge $e_0$ of the hexagonal lattice pointing to $v_0$
from outside of $G$.  Each coloring $A$ of the faces in $F$ determines an exploration tree
which is the union of exploration paths $T_v(A)$ over all vertices $v$ in $G$.

Let $f_0$ be the face in $F$ that is incident to $v_0$.  Now, given a subset $A$ of $F$, we define
a height function $h_A: F \rightarrow \mathbb Z$, defined up to additive constant, by writing
$h_A(f) - h_A(f_0)$ to be the number of left turns minus the number of right turns taken by the
sequence of edges in $T_v(A)$, where $v$ is the minimal vertex on $F$ in $T(A)$ (i.e., the first
vertex of $F$ hit by the exploration tree).   We refer to the value $\frac{2 \pi}{6} h_A(f)$ as the {\bf winding
number} of the face $f$.  Its value modulo $2\pi$ determines the angle of the edge of $T(A)$ that
points to the $\prec$-minimal vertex of $f$.

In Section \ref{s.openproblems}, we will make a few conjectures involving
these height function converges to the multiple of the Gaussian free field when $h_A(f)$.
We derive two simple combinatorial properties of the height function here.

\begin{proposition} If $f_1$ and $f_2$ are faces that share an edge, then $|h_A(f_1) - h_A(f_2)|
\leq 6$. \end{proposition}
\begin{proof}
Let $v$ be the minimal vertex of $f_1$ and $w$ the minimal vertex of $f_2$ in the $\prec$ ordering.
The proposition follows by examining the possible positions of $v$ and $w$.  \end{proof}

\begin{proposition} \label{monotonicityproposition} If $h_A(f_0)$ is fixed, then the height function $h_A$
is a decreasing function of $A$, i.e., $A \subset B$ implies $h_A \geq h_B$.  \end{proposition}

At first glance, the direction of the proposition may be surprising: it says that if we add faces
to $A$ (so that we {\em increase} the number of left turns at places where $T(A)$ hits faces for the
first time) then the height function goes down (i.e., we {\em decrease} number of left turns minus
right turns in $T(A)$ before the first time it hits a given face).  To get the right intuition,
consider the case that $G$ is the graph shown in Figure \ref{explorationtreefigure} and $A = F$; in
this case, by turning left each time we hit a hexagon for the first time, we force $T(A)$ to be a
single path that hugs the outer boundary and spirals clockwise inward from $v_0$. When $T(A)$ hits
a hexagon for the first time, the winding number is at its lowest possible value.

\begin{proof}
It is enough to prove that $h_A(f) \geq h_B(f)$ when $B \backslash A$ consists of a single face
$f'$.  Clearly, $T(A)$ and $T(B)$ will have the same minimal vertex $v'$ on $f'$, so we may as well
remove the path from $v_0$ to $v'$ and assume without loss of generality that $f' = f_0$, in which
case the result is a consequence of Lemma \ref{rotationlemma} (below).  \end{proof}

Sometimes it will be useful to compare $h_A$ with the height function $h_A'$ corresponding to a
modified exploration tree that begins at a boundary vertex $v_0' \not = v_0$ instead of at $v_0$.  To
make this comparison, we may choose the additive constant $h_A(f_0')$ to be $h_A(f_0)$ plus the
number of left turns minus the number of right turns in a path from the $e_0$ pointing to $v$ to
the $e_0'$ pointing to $v_0$ that doesn't intersect the interior of $G$ and makes a partial
revolution around $G$ in the clockwise direction.  We refer to such a change as a {\bf partial
clockwise rotation} of the model. We define counterclockwise rotations analogously.

\begin{lemma} \label{rotationlemma} If $A$ is fixed, and we replace $v_0$ with another vertex
$v_0'$ on the boundary and define $h_A'$ accordingly via a counterclockwise rotation, then $h_A'
\geq h_A$. \end{lemma}

\begin{proof}  We aim to prove $h_A'(f) \geq h_A(f)$ for some given face $f \in F$.  To this end,
let $C$ be the largest cluster of either white or black faces which
contains $f$ in its interior. (For example, if $f$ is the lone black
face surrounded by white faces in the center of Figure
\ref{explorationtreefigure}, then $C$ is large white cluster that
has five hexagons on the boundary of $G$.)  Write $\partial C$ for
the set of faces of $C$ which are incident to a hexagon in the
unbounded component of $\mathcal H \backslash C$.  If $f \not \in
\partial C$, then let $\tilde C$ be the component of $F \backslash
C$ containing $f$.  Let $P$ be the path in $T(A)$ connecting $v_0$
to the first vertex at which $T(A)$ hits $f$.  If $v_0$ does not lie
on $C$, then let $D$ be the cluster---with color opposite to that of
$C$---which borders $C$ and lies in the component of $F \backslash
C$ along which $v_0$ lies.  It is not hard to see that after $P$
hits $D$, it will follow the outer contour of $D$ counterclockwise
(if $D$ is black) or clockwise (if $D$ is white) until it first hits
$C$; thus, the first time $P$ hits $C$ will be at the first boundary
vertex of $C$ that lies counterclockwise (if $D$ is black) or
clockwise (if $D$ is white) of $v_0$ along the boundary of $G$.

Either way, once $P$ hits that vertex, it traces the complete contour of $C$ once (in a direction
depending on the color of $C$)---hitting all but one edge of the contour---before turning to hit a
face in $\tilde C$ (if $f \not \in \partial C$).  Since the position at which $P$ first hits $C$
(measured either clockwise or counterclockwise around the boundary of $C$) is monotone in the
position of $v_0$ (measured either clockwise or counterclockwise around the boundary of $G$), we
lose no generality in assuming that $C = G$, and if $f \in \partial C$, the result easily follows.
If $f \not \in \partial C$, then the position at which $P$ first hits $\tilde C$ (measured either
clockwise or counterclockwise around the boundary of $\tilde C$) is also monotone in the position
of $v_0$, so it is now enough to prove the result for $G = \tilde C$.  The lemma follows by
induction on the size of $C$.  \end{proof}

Proposition \ref{monotonicityproposition} implies that if we sample $A$ according to any measure
that satisfies the FKG inequality (i.e., increasing functions of $A$ are not negatively
correlated), then the random height function $h_A$ also satisfies the FKG inequality. In
particular, this is the case if we sample $A$ using Bernoulli percolation or a ferromagnetic Ising model.

\begin{figure}[ht]\label{bondpercolationfigure} \epsfbox[-30 22 295 155]{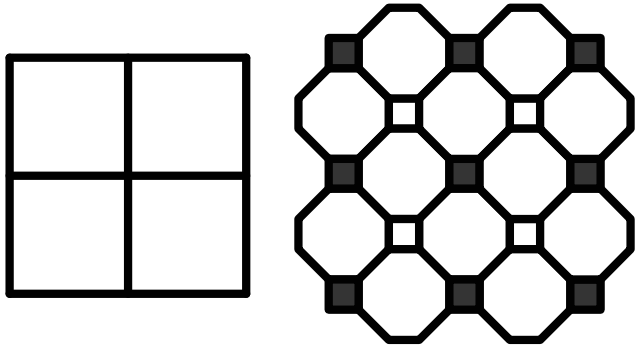} \caption {A grid
graph with nine vertices and twelve edges for bond percolation
(left) and a corresponding three-regular graph (right).  Each
percolation configuration of the left graph (i.e., a subset of the
edges) corresponds to a coloring of the octagons (one for each edge)
in the right graph.  The squares are given a fixed coloring in the
pattern shown.} \end{figure}

\subsection{Other lattices and bond percolation} \label{otherlattices}

The construction of the exploration tree and height functions
in the previous sections works for any periodic, three-regular planar,
periodic graph.  If the graph is not three-regular, we can make it
three-regular by replacing each degree $d > 3$ vertex with a
$d$-gon, and then coloring all of these $d$-gons white (i.e., the
$d$-gons are not allowed to be subsets of $A$, since they were not
present in our original graph).

In fact, we can also define loop ensembles and exploration trees
corresponding to instances of bond percolation
on a periodic lattice $G$. First form a three-regular graph $G'$
whose faces correspond to the edges, vertices, and faces of $G$;
then deterministically color those faces corresponding to vertices of $G$
black and those faces corresponding to faces of $G$ white.  Figure
\ref{bondpercolationfigure} illustrates this construction for a
small grid graph.

\section{Open problems} \label{s.openproblems}

The following two questions might be closely related; it is not
clear which of the two will be easier to address first.  (Recall
Section \ref{s.onmodel} for $O(n)$ model definitions.) When presenting this and other
conjectures about scaling limits, we will not specify the desired topology
of convergence, since part of the problem may be determining which
topology is most natural and tractable.

\begin{problem}
Are the scaling limits of the $O(n)$ models actually given by
\CLEk/, where $\kappa$ is as given in Section \ref{s.onmodel}?
Do the height functions of these models have scaling limits
given by a multiple of the Gaussian free field, with some boundary
conditions?  We conjecture that the answer to the second question
is yes whenever $\kappa > 4$ and $\beta = 1$.
\end{problem}

\medskip

\medskip

Next, suppose that $8/3 < \kappa \leq 4$.  Given an instance $\L$ of \CLEk/, we may
choose an orientation for each loop. We define a function $h_k(z)$ to be
the number of loops in the set $\{L^z_1, \ldots L^z_k\}$ that are
oriented counterclockwise around $z$ minus the number that are oriented
clockwise.

\medskip

\begin{problem}
What can be said about $h = \lim_{i \to \infty} h_i$?  We expect this
convergence to hold in the space of distributions, i.e.,
the limit of $\int_D \phi(z) h_i(z) dz$ should
exist almost surely for each smooth function $\phi$ on $D$.  When $\kappa = 4$,
the random distribution $h$ should be a multiple of the Gaussian free field.
What about the other values of $\kappa$?  Is there a natural description of these
fields that does not involve \SLE/?  Can we make sense of the
expectation of $\int_D \phi(z) h(z) dz$ when $h$ has piecewise constant boundary conditions?
Is the set of loops completely determined by the distribution $h$ almost surely?
\end{problem}

\medskip

In principle, our construction of \CLEk/ in terms of branching \SLEkr/ should allow one to compute
multi-point correlation functions
for the fields mentioned above, but it is not clear whether this can be done explicitly (or
how enlightening the answer will be).  Next, one may try to generalize Conjecture \ref{c.reversible}
as follows:

\medskip

\begin{problem}
Fix a domain $D$ with boundary points $a$ and $b$.
For what values of $\kappa$, $\rho$, $\mu$, and skew constant $\beta$ is it
the case that \SLEbkr/ is a continuous path almost surely?  When is
the law of an \SLEbkr/ from $a$ to $b$ in $D$ the same (up to parameterization) as the law
of its image under an anti-conformal map of $D$ that maps $b$ to $a$ and $a$ to $b$?
\end{problem}

\medskip
As discussed in Section \ref{s.confsymmetry}, a proof of Conjecture \ref{c.reversible} would
in particular yield a proof that the \CLEk/ loops are almost surely continuous when $D$
is a Jordan domain and $4 < \kappa < 8$.  One might expect that even in more general domains,
the \CLEk/ are almost surely continuous for these values of $\kappa$.  Some loops intersect the
boundary of $D$ almost surely, but it may be the case that the ``bad'' boundary points of
$D$ are rare enough that the loops are unlikely to intersect $\partial D$
at those boundary points.

\medskip
\begin{problem}
Let $D$ be an arbitrary simply connected planar domain.  Are all of the loops in a \CLEk/ almost
surely continuous in this case when $4 < \kappa < 8$?
\end{problem}

\medskip

If a path $\gamma$ chosen from \SLEr \kappa-6/ is almost surely
continuous, and $\kappa \leq 4$, then it is natural to define the
{\bf trunk} of $\gamma$ by $\{\gamma(t): O_t = W_t \}$.  Our
intuitive picture of \SLEkr/ is that it consists of the trunk
together with a pairwise disjoint collection of loops of the \CLE/,
each of which is rooted at a single point on the trunk.  If the skew
constant $\beta$ is $1$, then we expect all of the loops to lie to
one side of the trunk.  Otherwise, we expect there to be loops on
both sides of the trunk, where the fraction of loops which lie on
one side or the other is determined by $\beta$.

\medskip

\begin{problem}
For what values of $\kappa$, $\rho$, and skew constant $\beta$ is it
the case that the trunk of an \SLEr \kappa-6/ is almost surely
continuous?  Is the trunk also an \SLEkr/ process?  We conjecture that
when $\beta = 0$, the trunk has the law of an
\SLEkk \kappa', \frac{\kappa'-6}{2}, \frac{\kappa'-6}{2}/ process, where
$\kappa' = 16/\kappa$.
\end{problem}

\begin{figure}[t]\label{f.discretegasket}
%\epsfbox[-45 43 100 240]
\begin{center}
\rotatebox{90}{\reflectbox{\scalebox{.5}{\includegraphics{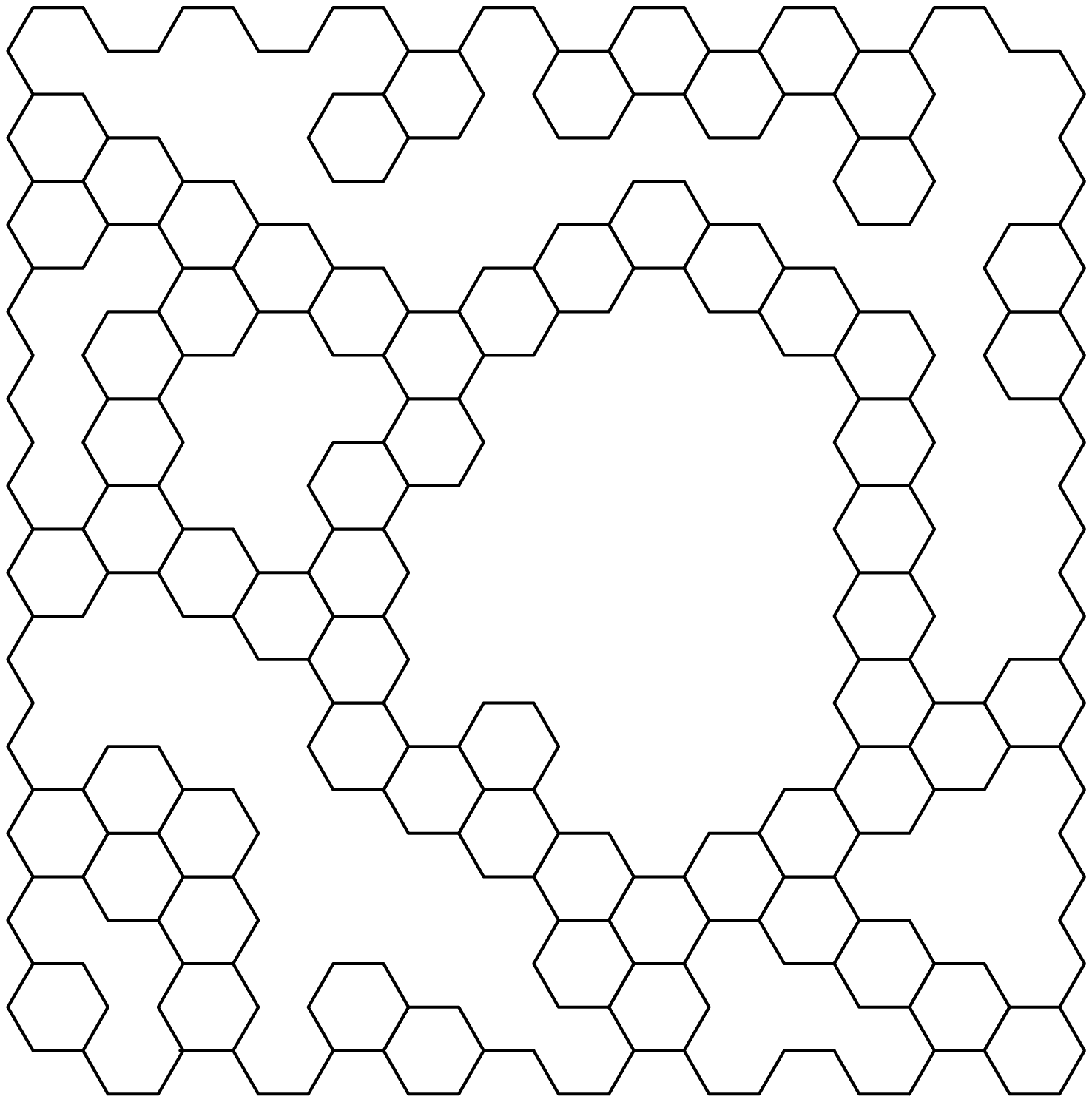}}}}
\caption {The discrete gasket of the loop ensemble in Figure \ref{explorationtreefigure}.}
\end{center}
\end{figure}

Given a disjoint simple loop ensemble in a planar graph, we define the
{\bf discrete gasket} to be the graph
obtained by deleting every vertex (plus incident edges) that is surrounded by a loop.
See Figure \ref{f.discretegasket}.  All of the vertices in the discrete gasket have degree two
or three---hence if we consider a coloring of the faces of the discrete gasket, we can draw
an exploration tree corresponding to that coloring.  A path in this exploration tree has
``faces'' on either side of it that correspond to loops of the $O(n)$ model.  The trunk of
an \SLEbr \kappa-6/, with $\mu = \beta = 0$, has \CLEk/ loops on either side of it.  This
suggests the following question:

\begin{problem} Consider independent Bernoulli percolation with $p=1/2$ on the faces of a gasket derived from a critical
$O(n)$ model (whose scaling limit we expect to be \CLEk/ for some $8/3 < \kappa \leq 4$).  What
is the scaling limit of the set of cluster boundaries of this percolation?  What is the scaling
limit of a branch of the exploration tree (say, from one fixed boundary vertex $a$ to another
boundary vertex $b$)?  We conjecture that its law is
the same as that of the trunk of a branching \SLEbr \kappa-6/ process with $\beta = \mu = 0$.
\end{problem}

\medskip

A modification of the above is that we independently color each face
comprised of multiple hexagons black with probability $p$ and independently color each face made
up of exactly one hexagon black with some probability $p'$.

\medskip

\begin{problem} For
each given value of $p$, is there a unique $p'$ for which the scaling limit of the exploration tree of
this coloring is given by the trunk of branching \SLEbr \kappa-6/ with $\beta = 2p-1$?  What
can be said about Ising and $O(n)$ models on the discrete gasket?
\end{problem}

\medskip

The intuition is that the $p'$ weight of the small hexagons generates the L\'evy compensation which
was necessary to make \SLEbr \kappa-6/ well defined when $\beta \not = 0$.

\medskip

The {\bf FK cluster model} corresponding to an expansion of
the $q$-state Potts models, may
be viewed as a random subset of the edges of a planar graph (i.e., a random
instance of non-independent bond percolation; see, e.g., \cite{\NienhuisKagerSurvey}
for details), where the probability of a set of edges
is proportional to $(e^\beta - 1)^b q^c$ where $\beta$ is some constant, $b$ is the
number of edges, and $c$ is the number of connected components of the subgraph
of the original graph containing those edges (and all of the original vertices).
As discussed in Section \ref{otherlattices}, each such subset determines
a collection of loops, and it is commonly conjectured that for a critical value
$\beta_c$ this set of loops has a non-trivial scaling limit.
Following an analogous conjecture for \SLEk/ given in \cite{\NienhuisKagerSurvey},
we ask the following:

\begin{problem}
When $0 < q \leq 4$ and $\beta = \beta_c$, is the scaling limit of the set of loops corresponding
to the critical FK clusters on a planar lattice given by \CLEk/ where $q = 2 + 2 \cos(8 \pi/ \kappa)$,
and $4 \leq \kappa \leq 8$?  Is this true in both the case of free boundary conditions (where
all subgraphs are allowed) and wired boundary conditions (where all boundary edges are deterministically
included in each subgraph---but the outer boundary of the outermost cluster is not counted as a loop)?
(The two are equivalent for self-dual graphs like $\mathbb Z^2$. \cite{\NienhuisKagerSurvey})
\end{problem}

Let $\L_j$ be the set of $j$-th nested loops in an instance of \CLEk/ (i.e., the set of loops of the
form $L^z_j$).  In the case of free boundary conditions, we may define a {\bf continuum FK cluster} $C$ to be
the set of points on or surrounded by a loop $L \in \L_j$, where $j$ is odd, minus the set of
all points surrounded by loops of the form $L^z_{j+1}$.  Each continuum FK cluster is a random
closed set.  In the case wired boundary conditions, a continuum FK cluster is the set of points on or
surrounded by a loop $L \in \L_j$, where $j$ is even (where we formally define $\L_0$ to
consist of the single loop given by the boundary of the domain), minus the set of
all points surrounded by loops of the form $L^z_{j+1}$.

In the discrete setting, one way to sample from the $q$-state Potts model is to first sample a collection
of FK clusters according to the model described above and then assign one of
the $q$ spins (uniformly at random and independently) to each cluster (assigning all
of the vertices in that cluster the corresponding spin).  (Free and wired boundary conditions
in the FK cluster model corresponding to free and constant-spin boundary conditions
in the corresponding Potts model.)  We now seek to define a continuum analog
of this construction.  In the continuum setting, we can also
uniformly and independently assign one of the $q$-states to each continuum FK cluster.  We then define a
{\bf continuum spin cluster} to be a connected component of the set of continuum FK clusters
of a given spin (with two continuum FK clusters considered adjacent if their intersection
is non-empty).

\medskip

\begin{problem} \label{prob.potts} We conjecture that the macrosopic same-spin clusters in the $q$-state Potts models
for $q \in \{2,3,4\}$ have scaling limits given by the continuum spin clusters described above.
\end{problem}

\medskip

Even for non-integer $1 < q \leq 4$, we can define the ``outermost spin cluster''
in the wired case to be the cluster of FK clusters consisting of those whose spins
are the ``same as the outermost cluster,'' where each cluster is assigned to have the same
spin as the outermost cluster with probability $1/q$.  In addition to discrete questions
like Problem \ref{prob.potts}, we can now ask a purely continuum question.

\begin{problem}
In the case of wired boundary conditions, is the law of the outermost continuum spin cluster corresponding
to \CLEk/ (for $\kappa \in [4,6)$, $q \in (1,4]$) given by the \CLEkk \kappa'/ gasket for $\kappa' = 16/\kappa$?
\end{problem}

\bibliographystyle{abbrv}
\bibliography{cle, mr, notmr}

\begin{thebibliography}{10}

\bibitem{MR1712629}
M.~Aizenman and A.~Burchard.
\newblock H\"older regularity and dimension bounds for random curves.
\newblock {\em Duke Math. J.}, 99(3):419--453, 1999.

\bibitem{math.Pr/0504036}
F.~Camia and C.~M. Newman.
\newblock {The Full Scaling Limit of Two-Dimensional Critical Percolation},
  2005.

\bibitem{CZ}
J.~Cardy and R.~M. Ziff.
\newblock Exact results for the universal area distribution of clusters in pe
  rcolation, {Ising} and {Potts} models, 2002.
\newblock arXiv:cond-mat/0205404.

\bibitem{MR2190302}
M.~Decamps, M.~Goovaerts, and W.~Schoutens.
\newblock Asymmetric skew {B}essel processes and their applications to finance.
\newblock {\em J. Comput. Appl. Math.}, 186(1):130--147, 2006.

\bibitem{D:vesicle}
B.~Duplantier.
\newblock Exact fractal area of two-dimensional vesicles.
\newblock {\em Physical Review Letters}, 64(4):493, 1990.

\bibitem{MR2065722}
W.~Kager and B.~Nienhuis.
\newblock A guide to stochastic {L}\"owner evolution and its applications.
\newblock {\em J. Statist. Phys.}, 115(5-6):1149--1229, 2004.

\bibitem{MR2129588}
G.~F. Lawler.
\newblock {\em Conformally invariant processes in the plane}, volume 114 of
  {\em Mathematical Surveys and Monographs}.
\newblock American Mathematical Society, Providence, RI, 2005.

\bibitem{MR2044671}
G.~F. Lawler, O.~Schramm, and W.~Werner.
\newblock Conformal invariance of planar loop-erased random walks and uniform
  spanning trees.
\newblock {\em Ann. Probab.}, 32(1B):939--995, 2004.

\bibitem{MR1466546}
J.~Pitman.
\newblock Partition structures derived from {B}rownian motion and stable
  subordinators.
\newblock {\em Bernoulli}, 3(1):79--96, 1997.

\bibitem{MR2000h:60050}
D.~Revuz and M.~Yor.
\newblock {\em Continuous martingales and {B}rownian motion}, volume 293 of
  {\em Grundlehren der Mathematischen Wissenschaften [Fundamental Principles of
  Mathematical Sciences]}.
\newblock Springer-Verlag, Berlin, third edition, 1999.

\bibitem{MR2153402}
S.~Rohde and O.~Schramm.
\newblock Basic properties of {SLE}.
\newblock {\em Ann. of Math. (2)}, 161(2):883--924, 2005.

\bibitem{math.Pr/0505368}
O.~Schramm and D.~B. Wilson.
\newblock {SLE coordinate changes}.
\newblock {\em New York J. Math.}, 11:659--669, 2005.

\bibitem{SchrammSheffieldWilson}
O.~Schramm; S.~Sheffield and D.~Wilson.
\newblock {CLE gaskets and conformal radii of CLE loops}.

\bibitem{SheffieldWerner}
S.~Sheffield and W.~Werner.
\newblock {Loop soup clusters and simple CLEs}.

\bibitem{MR1851632}
S.~Smirnov.
\newblock Critical percolation in the plane: conformal invariance, {C}ardy's
  formula, scaling limits.
\newblock {\em C. R. Acad. Sci. Paris S\'er. I Math.}, 333(3):239--244, 2001.

\bibitem{MR1335470}
S.~Watanabe.
\newblock Generalized arc-sine laws for one-dimensional diffusion processes and
  random walks.
\newblock In {\em Stochastic analysis (Ithaca, NY, 1993)}, volume~57 of {\em
  Proc. Sympos. Pure Math.}, pages 157--172. Amer. Math. Soc., Providence, RI,
  1995.

\bibitem{MR2023758}
W.~Werner.
\newblock S{LE}s as boundaries of clusters of {B}rownian loops.
\newblock {\em C. R. Math. Acad. Sci. Paris}, 337(7):481--486, 2003.

\bibitem{MR2079672}
W.~Werner.
\newblock Random planar curves and {S}chramm-{L}oewner evolutions.
\newblock In {\em Lectures on probability theory and statistics}, volume 1840
  of {\em Lecture Notes in Math.}, pages 107--195. Springer, Berlin, 2004.

\end{thebibliography}

%\bibliographystyle{halpha}
%\addcontentsline{toc}{section}{Bibliography}
%\bibliography{mr,prep,notmr}

\bigskip

\filbreak
\begingroup
\small
\parindent=0pt

\bigskip
\vtop{
\hsize=1.3in
Courant Institute\\
251 Mercer Street\\
New York, NY 10012 \\
\\{\tt {sheff@math.nyu.edu}}}
\endgroup
\filbreak
\end{document}